\documentstyle{amsppt}
\def\today
{\ifcase\month\or
     January\or February\or March\or April\or May\or June\or
     July\or August\or September\or October\or November\or December\fi
     \space\number\day, \number\year}
%\UseAMSsymbols
\hsize 5.5 true in
\vsize 8.5 true in
\parskip=\medskipamount
\NoBlackBoxes

\def\mathbb{\Bbb}

\def\mathcal{\Cal}

\def\mod{\text{\rm mod\,}}

\def\ve{\varepsilon}
\def\vp{\varphi}
\def\arrowk{^\to{\kern -6pt\topsmash k}}
\def\arrowK{^{^\to}{\kern -9pt\topsmash K}}
\def\arrowr{^\to{\kern-6pt\topsmash r}}
\def\bark{\bar{\kern-0pt\topsmash k}}
\def\arrowvp{^\to{\kern -8pt\topsmash\vp}}
\def\arrowf{^{^\to}{\kern -8pt f}}
\def\arrowg{^{^\to}{\kern -8pt g}}
\def\arrowu{^{^\to}a{\kern-8pt u}}
\def\arrowt{^{^\to}{\kern -6pt t}}
\def\arrowe{^{^\to}{\kern -6pt e}}
\def\tk{\tilde{\kern 1 pt\topsmash k}}
\def\barm{\bar{\kern-.2pt\bar m}}
\def\barN{\bar{\kern-1pt\bar N}}
\def\barA{\, \bar{\kern-3pt \bar A}}

\def\mathbb{\Bbb}

\NoRunningHeads
\TagsOnRight
\topmatter
\title
On Representation of Integers by Binary Quadratic Forms
\endtitle
\author
J.~Bourgain, E.~Fuchs
\endauthor

\abstract
Given a negative $D>-(\log X)^{\log 2-\delta}$, we give a new upper bound on the number of square free integers $<X$ which are represented by some but not all forms of the genus of a primitive positive definite binary quadratic form $f$ of discriminant $D$.  We also give an analogous upper bound for square free integers of the form $q+a<X$ where $q$ is prime and $a\in\mathbb Z$ is fixed.  Combined with the $1/2$-dimensional sieve of Iwaniec, this yields a lower bound on the number of such integers $q+a<X$ represented by a binary quadratic form of discriminant $D$, where $D$ is allowed to grow with $X$ as above.  An immediate consequence of this, coming from recent work of the authors in [BF], is a lower bound on the number of primes which come up as curvatures in a given primitive integer Apollonian circle packing.
\endabstract
\address
Institute for Advanced Study,
1 Einstein Drive, Princeton, NJ 08540
\endaddress
\address
Institute for Advanced Study,
1 Einstein Drive, Princeton, NJ 08540
\endaddress
\thanks The first author was partially supported by NSF Grants DMS-0808042 and DMS-0835373
\endthanks
\thanks The second author was supported by NSF Grant DMS-0635607
\endthanks
\endtopmatter

\document

\noindent
{\bf \S0. Introduction}

Let $f(x, y)=ax^2+bxy+cy^2\in\Bbb Z[x, y]$ be a primitive positive-definite binary quadratic form of negative discriminant $D=b^2-4ac$.
For $X\to \infty$, we denote by $U_f(X)$ the number of positive integers at most $X$ that are representable by $f$.
The problem of understanding the behavior of $U_f(X)$ when $D$ is not fixed, i.e. $|D|$ may grow with $X$, has been addressed in several recent papers, in
particular in [Bl] and [B-G].
What is shown inthese papers, on a crude level, is that there are basically three ranges \big(we restrict ourselves to discriminants satisfying $\log |D|\leq \text{O} (\log\log
X$)\big)
\medskip
\roster
\item
"{(i)}"  $|D| \ll (\log X)^{(\log 2)-\ve}$. Then $U_f (X) \gg_\ve X(\log X)^{-\frac 12-\ve}\hfill (0.1)$
\medskip
\item
"{(ii)}"  $|D|\gg (\log X)^{2(\log 2)+\ve}$. Then  $U_f(X) \asymp \frac X{\sqrt D}\hfill (0.2)$
\medskip
\item
"{(iii)}"  The intermediate range.
\endroster

As Blomer and Granville explain in [B-G], this transitional behavior is due to the interplay between the size $h$ of the class group $\Cal C$ and the typical number of prime
factors of an integer $n\sim X$.
A precise elaboration of the underlying heuristics was kindly communicated by V.~Blomer to the authors and is reproduced next.
The number of integers $n<X$ with $k$ prime factors $p$ split in the quadratic number field (i.e. $(\frac Dp) =1$) is of the order
$$
\frac X{\log X} \ \frac 1{2^k} \ \frac {(\log\log X)^{k-1}}{(k-1)!}.\tag 0.3
$$
Note that summation of (0.3) over $k$ gives $\frac X{\sqrt {\log X}}$, which is the number of integers at most $X$ represented by some form of discriminant $D$.

Moreover, applying Stirling's formula, we see that the main contribution comes from integers with $k\sim\frac 12\log\log X$ prime factors.

Next, ignoring ambiguous classes, these $k$ primes yield $2^k$ classes (with possible repetition) in $\Cal C$ that represent the given integer $n$.
Hence, roughly speaking, one would expect that typically $n$ is represented by each class of its genus provided $2^k\gg h$, which amounts to
$$
h<(\log X)^{\frac {\log 2} 2 -\ve}\tag 0.4
$$
corresponding to alternative (i).

On the other hand, if $D$ is sufficiently large, the $2^k$ classes will be typically distinct.
Assuming some mild form of equidistribution in the class group when varying $n$, we expect for the number of integers $n < X$ with $k$ prime
factors represented by a given class to be of order
$$
\frac {2^k}h \cdot (0.3) =\frac X{h\log X} \ \frac {(\log\log X)^{k-1}}{(k-1)!}
\tag 0.5
$$
with total contribution $\text{O}\Big(\frac Xh\Big)$, attained when $k\sim\log\log X$ (at this level of the discussion, there is no difference between $h$ and
$\sqrt D$).

In this paper, we consider only the lower range (i).
Our aim is to substantiate further the heuristic discussed above according to which, typically, all classes of the genus of $n\sim X$,  $n$
representable by a form of discriminant $D$, do actually represent $n$.

More precisely, we prove the following (as consequence of Theorem 2 in [B-G]).

\proclaim
{Theorem 2'}
Let $D$ be a negative discriminant satisfying
$$
|D|<(\log X)^{\log 2-\delta}\tag 0.6
$$
for some fixed $\delta>0$.
Then there is $\delta'=\delta'(\delta)>0$ such that

\noindent
$\#\{n\sim X$; $n$ square free, representable by some form of discriminant $D$ but not by all forms of the genus\}
$$
< \frac X{(\log X)^{\frac 12+\delta'}}.\tag 0.7
$$
\endproclaim

Note that though [Bl], [B-G] establish (0.1) (in fact in a more precise form, cf. Theorem 5 in [B-G]), their results do not directly
pertain to the phenomenon expressed in Theorem 2'.
As pointed out in [B-G], it was shown on the other hand 
by Bernays that almost all integers represented by some form in a given genus can be represented by
all forms in the genus, but assuming the much stronger restriction
$$
D\ll (\log\log X)^{\frac 12-\ve}.\tag 0.8
$$
A result in the same spirit was also obtained by Golubeva [Go].

The proof of Theorem 2' rests on a general result from arithmetic combinatorics (Theorem 1 below) that we describe next.
Assume $G$ a finite abelian group $(G=\Cal C^2$ in our application) in which group operation will be denoted additively.
Given a subset $A\subset G$, we introduce the set
$$
s(A) =\Big\{\sum x_i; \{x_i\} \text{ are distinct elements of $A$}\Big\}. \tag 0.9
$$
The issue is then to understand what it means for $A$ that $s(A)\not = G$.
It turns out that there are basically two possibilities.
In the first, $A$ is contained, up to a bounded number of elements, in a proper subgroup $H$ of $G$ of bounded index [G:H].
The second scenario is as follows.
There are $k$ elements $x_1, \ldots, x_k\in A$ with
$$
k< (1+\ve)\frac {\log |G|}{\log 2}\tag 0.10
$$
and a subset $\Omega_{x_1, \ldots, x_k}\subset G$ (determined by $x_1, \ldots, x_k$), such that $A\subset\Omega_{x_1, \ldots, x_k}$ and
$$
|\Omega_{x_1, \ldots, x_k}|< \ve|G|\tag 0.11
$$
(we are assuming here that $|G|$ is large).

To prove Theorem 1, one applies the greedy algorithm.
Thus given $x_1, \ldots, x_k \in A$, we select $x_{k+1}\in A$ as to optimize the  size of $s(x_1, \ldots, x_{k+1})$.
If we do not reach $s(x_1, \ldots, x_k)=G$ with $k$ satisfying (0.10), then
$$
A\subset \{x_1, \ldots, x_k\}\cup \Omega\tag 0.12
$$
where the elements $x\in\Omega$ have the property that
$$
|s(x_1, \ldots, x_k, x)|\approx |s(x_1, \ldots, x_k)|.\tag 0.13
$$
Assuming $\Omega$ fails (0.11), the first alternative is shown to occur.
The argument involves combinatorial results, such as a version of the Balog-Szemeredi-Gowers theorem and also Kneser's theorem.
The reader is referred to the book [T-V] for background material on the matter.

Once Theorem 1 is established, deriving Theorem 2 is essentially routine.
We make use, of course, of Landau's result [L2] (established in [Bl] with uniformity in the discriminant), on the distribution of the primes represented
by a given class $C\in \Cal C$ -- namely, for $\Cal P_C$ the set of primes represented by a class $C$, 
$$
|\{p\in\Cal P_C; p\leq \xi\}|=\frac 1{\ve(C)h} \int_1^\xi \frac {dt}{\log t}+ C(\xi e^{-c\sqrt{\log \xi}})\tag 0.14
$$
for $\xi\to\infty$, with $\ve(C) =2$ if $C$ is ambiguous and $\ve(C)=1$ otherwise.

The nontrivial upper bound (0.7) is then obtained by excluding certain additional prime divisors, i.e. satisfying $\big(\frac Dp\big)\not= -1$, using
standard upper bound sieving.

The same approach permits to obtain a similar result considering now shifted primes, i.e. integers $n$ of the form $n=a+q$ with $a$ fixed and
$q$ a prime number. Thus

\proclaim
{Theorem 3'} Under the assumption (0.6), fixing $a\in\mathbb Z$, we have

\noindent
$|\{q+a\sim X$; $q$ prime, $q+a$ squarefree representable by some form of discriminant $D$ but not by all forms of the genus\} $|$ 
$$
<\frac X{(\log X)^{\frac 32+\delta'}}.\tag 0.15
$$
\endproclaim

On the technical side, only crude sieving bounds are needed for our purpose and they can be obtained by the simple inclusion/exclusion principle without
the need of Brun's theory.  The arguments covering the specific problem at hand were included in the paper (see Lemmas 4 and 5), which turned out to be
more convenient than searching for a reference.
Note that the proof of Lemma 5 involves sieving in the ideals and the required remainder estimates are provided by Landau's extension of the
Polya-Vinogradov inequality for Hecke characters [L1].

The motivation behind Theorem 3' lies in a result due to H.~Iwaniec [Iw] on the number of shifted primes that are representable by the genus of a quadratic form.  This in turn is applicable to counting primes which appear as curvatures in a primitive integer Apollonian circle packing using a method similar to that in \cite{BF}, where the authors prove that the integers appearing as curvatures in a primitive integer Apollonian packing make up a positive fraction of $\mathbb Z$.

Specifically, let $P$ be a primitive integer Apollonian packing, and let $a\not=0$ denote a curvature of a circle in $P$.  From \cite{BF}, we have that the set $S_a$ of integers less than $X$ represented by certain shifted binary quadratic forms $f_a(x,y)-a$, where the discriminant $D(f_a)=-4a^2$, correspond to curvatures of circles in $P$.  Let $\frak P_a\subset S_a$ denote the set of primes in $S_a$.  We may then compute a lower bound for the number of primes less than $X$ appearing as curvatures in $P$ by bounding
$$\big|\bigcup_{a}\frak P_a\big|$$
where the $a$'s range over a set of our choice.  The aim is to use the $\frac{1}{2}$-dimensional sieve of Iwaniec to first determine the cardinality of $\frak P_a$.  In [Iw], Iwaniec proves upper and lower bounds for the number of primes less than $N$ represented by $\phi(x,y)+A$, where $\phi(x,y)$ is a positive definite binary quadratic form and $A$ is an integer.  He shows
$$\frac{X}{(\log X)^{3/2}}\ll S(X,\phi,a)\ll \frac{X}{(\log X)^{3/2}}$$
where $S(X,\phi,a)$ denotes the number of primes less than $X$ represented by $\phi(x,y)+A$.  Here the discriminant of $\phi$ is fixed, and the bounds above are obtained by considering the count over all forms in the genus of $\phi$: namely, for fixed discriminant, bounds for $S(X,\phi,a)$ are easily derived from bounds for
$$
S_1(X,\phi,a)=\sum_{\Sb{p\leq X}\\{(x,y)=1,f\in R_{\phi}}\\
{p=f(x,y)+a}\endSb}1
$$
where $R_{\phi}$ denotes the genus of $\phi$.  In order to apply this to finding bounds for $|\frak P_a|$ where $a$ is allowed to grow with $X$, we must understand both how $S_1(X,\phi,a)$ depends on the discriminant of $\phi$, and how $S$ relates to $S_1$ in the case that $D$ is not fixed.  The latter is explained by Theorem 3' for $D$ satisfying (0.6), while the former is done via a careful analysis of the dependence on the discriminant in [Iw] for $D<\log X$.  This is discussed in the Appendix.  Note that in the application to Apollonian packings, the discriminant of $\phi$ is always of the form $-4a^2$, but our results apply to a much more general discriminant.

Indeed, Theorem 1 in [Iw] combined with Theorem 3' above implies the following

\proclaim
{Corollary 4} Let $D<0$ satisfy (0.6) and $f$ be a primitive positive definite binary form of discriminant $D$.
Then
$$
|\{q+a\sim X; q \text { prime, } q+a \text { representable by $f$}\}|\gg \frac X{(\log X)^{\frac 32+\ve}}\tag 0.16
$$
\endproclaim
(we assume here $a\in\Bbb Z$ fixed for simplicity).

\noindent
{\bf Acknowledgement:} The authors are grateful to V.~Blomer for several private communications.

\bigskip

\noindent
{\bf \S1. A result in combinatorial group theory}

The aim of this section is to prove Theorem 1 below.
 
\noindent
{\bf (1).} Let $G$ be an abelian group, $|G|=h'$.

For $A\subset G$, denote

$$
s(A)=\Big\{\sum x_i; \text{ $\{x_i\}$ distinct elements of $A$}\Big\}\tag 1.0
$$
(the set of sums of distinct elements of $A$.)

Assume $s(A)\not= G$.
We would like to specify the structure of such $A$.

Start with the following algorithm.

Take $x_1\in A$.

Assume obtained $x_1, \ldots, x_j$, take $x_{j+1}$ as to maximize
$$
s(x_1, \ldots, x_{j+1}).
$$
Let $\delta_j= \frac {|s(x_1, \ldots, x_j)|}{h'}$.

One has
$$
\spreadlines{6pt}
\align
|s(x_1, \ldots, x_j, x)|&=|s(x_1, \ldots, x_j)\cup\big(s(x_1, \ldots, x_j)+x\big)|\\
&=2 \big|s(x_1, \ldots, x_j)\big|-\big|s(x_1, \ldots, x_j)\cap \big(s(x_1, \ldots, x_j)+x\big)\big|.
\endalign
$$
Hence
$$
\Bbb E_x[|s(x_1, \ldots, x_j, x)|]= 2\delta_j h'-\delta_j^2 h'=\delta_j(2-\delta_j)h'.\tag 1.1
$$

On the other hand, for all $x$
$$
|s(x_1, \ldots, x_j, x)|\leq (2\delta_j) h'.\tag 1.2
$$
Fix $\ve>0$.

For $\delta_j<\frac 12$, define
$$
\Omega=\{x\in G; |s(x_1, \ldots, x_j, x)|< (2-\ve)\delta_j h'\}.
$$
Then, from (1.1), (1.2)
$$
\Bbb E_x [|s(x_1, \ldots, x_j, x)|] \leq (2\delta_j h')\Big(1- \frac{|\Omega|}{h'}\Big)
+(2-\ve)\delta_j h'\frac{|\Omega|}{h'}
$$
implying
$$
|\Omega|<\frac {\delta_j}\ve h'.\tag 1.3
$$
For $\delta_j>\frac 12$, define
$$
\Omega=\{x\in G; |s(x_1, \ldots, x_j, x)|< \big(1-(1-\delta_j)^{3/2}\big)h'\}.
$$
Similarly
$$
|\Omega| < (1-\delta_j)^{1/2} h'.\tag 1.4
$$

It follows from (1.3) that either

(1.5) there exist $x_1, \ldots , x_k \in A$ s.t.
$$
|s(x_1, \ldots, x_k)|> \ve^2 h'\tag 1.6
$$
with
$$
k<\frac {\log h'}{\log 2-\frac \ve 2}\tag 1.7
$$
or

\noindent
(1.8)  \ There exists elements $x_1, \ldots, x_k \in A$ and a set
$ \Omega_{x_1, \ldots, x_k} \subset G$ satisfying
$$
\qquad A\subset \{x_1, \ldots, x_k\}\cup \Omega_{x_1, \ldots, x_k}\tag 1.9
$$
$$
\qquad k< \frac {\log h'}{\log 2 -\frac\ve 2}\tag 1.10
$$
$$
\qquad |\Omega_{x_1, \ldots, x_k}|< \ve h'.\tag 1.11
$$

\bigskip
\noindent
{\bf (2).} Let $A_1\subset A$ s.t.
$$
\delta h' = s(A_1)> \ve^2 h'.\tag 2.1
$$
Fix $\ve_1> 0$ and define
$$
\Omega=\{ x\in G; |s(A_1\cup \{x\})| < (1-\ve_1)|s(A_1)|+\ve_1 h'\}.\tag 2.2
$$
If $(A\backslash A_1)\cap \Omega^c\not= \emptyset$, we add an element and increase the density from $\delta$ in (2.1) to $(1-\ve_1)\delta+\ve_1$.

Assume this process can be iterated $r$ times.

We obtain a set $A_1'$ such that $s(A'_1)$ has density at least $\delta'$ satisfying
$$
1-\delta'= (1-\ve_1)^r (1-\delta)
$$
and thus $|s(A_1')|> (1-\ve^2) h'$ for
$$
r\sim\frac {\log \frac 1\ve}{\ve_1}.\tag 2.3
$$
Continuing the process with $A_1'$ and $\delta> 1-\ve^2>\frac 12$ gives a subset $A_1''\subset A$ so that $s(A_1'')= G$ and
$$
|A_1''|\leq \frac { \log h'}{(\log 2)-\ve} + c\frac {\log \frac 1\ve}{\ve_1} +\log\log h'\tag 2.4
$$
unless we are in alternative (1.8) with (1.10) replaced by (2.4).

 Thus it remains to analyze the case when the iteration fails.

If $|\Omega|< \ve h'$, we are again in the situation (1.8) with (1.10) replaced by 
$$
\frac {\log h'}{\log 2-\ve}+c\frac{\log\frac
1\ve}{\ve_1}.
$$

Assume next $\Omega$ defined in (2.2) satisfies
$$
|\Omega|> \ve h'.\tag 2.5
$$
Denoting $B= s(A_1')$, we have by (2.1) and definition of $\Omega$ that
$$
|B|> \ve^2 h'\tag 2.6
$$
and
$$
|B\cap (B+x)|>(1-\ve_1\ve^{-2})|B| \text { for } x\in\Omega.\tag 2.7
$$
Hence
$$
1_B * 1_{-B}>(1-\ve_1\ve^{-2})|B| \text { on }\Omega\tag 2.8
$$
implying in particular that
$$
|B|>(1-\ve_1\ve^{-2})|\Omega|.\tag 2.9
$$
\bigskip

\noindent
{\bf (3).} Assume (2.6)-(2.9).

Thus
$$
\langle 1_B, 1_B* 1_\Omega\rangle =\langle 1_B * 1_{-B}, 1_\Omega\rangle \geq (1-\ve_1 \ve^{-2}) |B|\, |\Omega|
\tag 3.1
$$
and
$$
\Vert 1_B*1_\Omega\Vert_2\geq (1-\ve_1\ve^{-2})|B|^{\frac 12}|\Omega|.
$$
Squaring and using the fact that $\Omega$ is symmetric
$$
\Vert 1_B*1_\Omega *1_\Omega\Vert_2\geq (1-\ve_1\ve^{-2})^2 |B|^{\frac 12}|\Omega|^2
$$
and for any given $r$ (= power of 2)
$$
\Vert 1_B* 1_\Omega^{(r)}\Vert_2 \geq (1-\ve_1\ve^{-2})^r |B|^{\frac 12}|\Omega|^r.\tag 3.2
$$
(where $1_\Omega^{(r)}$ denotes the $r$ fold convolution).

We will rely on the following
\proclaim
{Lemma 1} Let $\mu$ be a probability measure on a discrete additive group $G$ and assume (for small $\kappa$)
$$
\Vert\mu * \mu\Vert_2 >(1-\kappa)\Vert\mu\Vert_2.\tag 3.3
$$
Then there is a subgroup $H$ of $G$ s.t.
$$
\frac 12\Vert\mu\Vert_2^{-2} <|H|< 2\Vert\mu\Vert_2^{-2}\tag 3.3'
$$
and for some $z\in G$
$$
\Big\Vert\mu - \frac {1_{H-z}}{|H|}\Big\Vert_1 < c\kappa ^{1/12}.\tag 3.3''
$$
\endproclaim

\noindent
{\bf Proof.} 

>From (3.3) we have
$$
\sum_x\Big|\sum_y \mu(x-y) \mu(y)\Big|^2>(1-\kappa)^2 \Vert\mu\Vert_2^2
$$
and
$$
\sum_{y_1, y_2} \langle \mu_{y_1}, \mu_{y_2}\rangle \mu(y_1)\mu(y_2) > (1-\kappa)^2 \Vert \mu\Vert_2^2
$$
implying
$$
\spreadlines{6pt}
\align
\sum \Vert \mu_{y_1} -\mu_{y_2}\Vert^2_2  \ \mu (y_1) \mu (y_2)&<  2\big(1-(1-\kappa)^2\big) \Vert\mu\Vert_2^2
\\
&< 4\kappa \Vert\mu\Vert_2^2.
\endalign
$$
Hence there is $y_0\in G$ s.t.
$$
\sum\Vert\mu_y-\mu_{y_0}\Vert^2_2 \ \mu(y)< 4\kappa\Vert\mu\Vert^2_2
$$
and by translation of $\mu$, we may assume $y_0 = 0$, thus
$$
\sum\Vert\mu_y -\mu\Vert^2_2 \ \mu(y) < 4\kappa \Vert\mu\Vert^2_2.
$$
Denote
$$
U =\{y\in G; \Vert\mu -\mu_y\Vert_2 < \kappa^{1/3} \Vert\mu\Vert_2\}.
$$
Hence, from the preceding
$$
\mu (G\backslash U)< 4\kappa^{1/3}.
$$
Since
$$
\frac 1{|U|}\sum_{y\in U} \Vert\mu -\mu_y\Vert_2 < 4\kappa^{1/3} \Vert\mu\Vert_2
$$
it follows by convexity that
$$
\Big\Vert \mu -\mu * \frac {1_U}{|U|} \Big\Vert_2 < 4\kappa^{1/3} \Vert\mu\Vert_2
$$
and in particular
$$
\Vert\mu\Vert_2 \leq\frac 1{|U|^{1/2}} + 4\kappa^{1/3} \Vert\mu\Vert_2
$$
$$
\Vert\mu\Vert_2< \frac {1+4\kappa^{1/3}}{|U|^{1/2}}.
$$
Next, write
$$
\spreadlines{6pt}
\align
\Big\Vert\mu -\frac {1_U}{|U|}\Big\Vert^2_2 &=\Vert\mu\Vert^2_2 +\frac 1{|U|} -2 \frac {\mu(U)}{|U|}\\
&\leq \frac {2+10\kappa^{1/3} -2 (1-4\kappa^{1/3})}{|U|}\\
&< \frac {18\kappa^{1/3}}{|U|}.
\endalign
$$
Hence
$$
\Big\Vert\mu -\frac {1_U}{|U|} \Big\Vert_2< \frac {5\kappa^{1/6}}{|U|^{1/2}}\tag 3.4
$$
and also
$$
\spreadlines{6pt}
\align
\Big\Vert\mu -\frac {1_U}{|U|}\Big\Vert_1& \leq \mu(U^c)+\sum_{x\in U} \Big|\mu(x)-\frac 1{|U|}\Big|\\
&\leq 4\kappa^{1/3}+|U|^{1/2} \Big\Vert\mu -\frac {1_U}{|U|}\Big\Vert_2\\
&< 6\kappa^{1/6}.\tag 3.5
\endalign
$$
>From (3.3), (3.4), (3.5), we have
$$
\Big\Vert\frac {1_U}{|U|} * \frac {1_U}{|U|} \Big\Vert_2 > (1-20\kappa^{1/6}) \frac 1{|U|^{1/2}}
$$
hence
$$
E_+(U, U)= \Vert 1_U*1_U\Vert^2_2 > (1-40 \kappa^{1/6}).|U|^3
$$
where $E_+$ refers to the additive energy.

We apply now some results from arithmetic combinatorics.

First, by (2.5.4), p.82 from [TV] (B-S-G in near-extreme case), there are subsets $U', U'' \subset U$ s.t.
$$
|U'|, |U''|>(1-10\kappa ^{1/12})|U|
$$
and
$$
|U'-U''|< (1+20\kappa^{1/12})|U|.
$$
Thus from Ruzsa's triangle inequality, also
$$
\spreadlines{6pt}
\align
|U'-U'|\leq \frac {|U'-U''|^2}{|U''|} &< (1+ 60\kappa^{1/12})|U|\\
&<(1+80\kappa^{1/12})|U|.\tag 3.6
\endalign
$$
Next, we apply Kneser's theorem (see [TV], Theorem 5.5, p. 200).

For $T\subset G$, denote
$$
\text{Sym}_1(T) =\{x\in G; T+x=T\}
$$
the symmetry group of $T$.

Then by Kneser's theorem, see [T-V]
$$
|T-T|\geq 2|T| -|\text {Sym}_1 (T-T)|
$$
and application with $T=U'$ gives
$$
|\text{Sym}_1 (U'-U')|> (1-80\kappa^{1/12})|U'|.\tag 3.7
$$
Denoting $H=\text { Sym}_1 (U'-U'), H\subset U'-U$ and thus
$$
\spreadlines{6pt}
\align
|H|\, |U'| &\leq\sum_{z\in U'-U'-U'}|H\cap (U'+z)|\\
&\leq |U'-U'-U'||\max_z |H\cap (U' +z)|\\
&< (1+300\kappa^{1/12})|U'|\max_z |H\cap (U'+z)|
\endalign
$$
from (3.6) and sumset inequalities. Therefore, there is some $z\in G$ s.t.
$$
|(H-z) \cap U'|> (1+300 \kappa^{1/12})^{-1} |H|
$$
and in view of (3.7)
$$
|U' \vartriangle (H-z)|< 1000 \kappa^{\frac 1{12}} |U|
$$
and
$$
|U\vartriangle (H-z)|< 1000 \kappa^{\frac 1{12}} |U|.\tag 3.8
$$
>From (3.5), (3.8) we have
$$
\Big\Vert\mu -\frac {1_{H-z}}{|H|}\Big\Vert_1 < C\kappa^{\frac 1{12}}.\tag 3.9
$$
>From (3.4), (3.8), we obtain (3.31) proving Lemma 1.
\bigskip

Returning to (2.5), (3.2), we have that
$$
\Big\Vert\Big(\frac {1_\Omega}{|\Omega|}\Big)^{(r)} \Big\Vert_2 
$$
decreases in $r$ and is between $\frac 1{\sqrt h'}$ and $\frac 1{\sqrt{\ve h'}}$.
Hence there is some $\tau$
$$
 \log r <\frac c\kappa \log\frac 1\ve\tag 3.10
$$
such that $\mu=(\frac{1_\Omega}{|\Omega|})^{(r)}$ satisfies (3.3).

>From (3.2), (3.3$'$), we conclude that
$$
\spreadlines{6pt}
\align
\Big\Vert 1_B*\frac {1_H}{|H|} \Big\Vert_2 &\geq \big((1-\ve_1 \ve^{-2})^r -c\kappa^{1/12})|B|^{1/2}\\
&> (1-c\kappa^{1/12})|B|^{1/2}\tag 3.11
\endalign
$$
provided
$$
\ve_1 <\Big(\frac 1\ve\Big)^{c\kappa^{-1}}.\tag 3.12
$$
Also, from (3.3') and the preceding
$$
|H|>  \frac 12|\Omega| >\frac\ve 2 h'.\tag 3.13
$$
Let $\{H_\alpha\}$ be the cosets of $H\subset G$. Then
$$
\Vert 1_B*1_H\Vert^2_2 =\sum_\alpha \Vert 1_{(B\cap H_\alpha)}*1_H\Vert^2_2.
$$
Let $\kappa_1>0$ be a small parameter and define
$$
I_0 =\{\alpha; |B\cap H_\alpha|> (1-\kappa_1) |H|\}
$$
and $I_1$ the complement.

One has
$$
\Vert 1_{(B\cap H_\alpha)}*1_H\Vert^2_2 = E_+(H, B\cap H_\alpha) \leq |B\cap H_\alpha|^2\cdot |H|
$$
and hence, by (3.11)
$$
\spreadlines{6pt}
\align
(1-c\kappa^{1/12})|B|\cdot |H|^2&\leq |H|\sum|B\cap H_\alpha|^2\\
&\leq |H|\Big(\sum_{\alpha\in I_0}|H| \, |B\cap H_\alpha|+ (1-\kappa_1) \sum_{\alpha\in I_1} |H| \, |B\cap H_\alpha|\Big)
\\
&\leq |H|^2 (|B|-\kappa_1\sum_{\alpha\in I_1} |B\cap H_\alpha|).
\endalign
$$
Hence $B=B_0\cup B_1$ with
$$
|B_1|=\sum_{\alpha\in I_1} |B\cap H_\alpha|< c\kappa^{1/12}\kappa_1^{-1}|B|.\tag 3.14
$$
Assume
$$
\kappa\ll \kappa_1^{12}\tag 3.15
$$
so that in particular $I_0\not= \emptyset$.

Let $y\in A\backslash A_1'$.
Then $y\in \Omega$ and by (2.7)
$$
|B\cap (B+y)|> (1-\ve_1\ve^{-2})|B|.
$$
Let $\vp:G\to G/H =I_0\cup I_1$.

If $\alpha\in I_0$, then
$$
\spreadlines{6pt}
\align
|\big((B\cap H_\alpha)+y\big) \cap B|&\geq |(B+y)\cap B| -\sum_{\alpha'\not=\alpha} |B\cap H_{\alpha'}|\\
&> (1-\ve_1\ve^{-2})|B|-|B|+|B\cap H_\alpha|\\
&> (1-\kappa_1) |H|-\ve_1 \ve^{-2}|B|\\
&\underset {(3.13)}\to > (1-\kappa_1 -2\ve_1\ve^{-3})|H|.
\endalign
$$
Thus certainly
$$
|H_{\alpha+\vp(y)} \cap B|> (1-\kappa_1-2\ve_1 \ve^{-3})|H|.
$$
>From (3.14), if $\beta\in I_1$
$$
|H_\beta\cap B|< c\kappa^{1/12}\kappa_1^{-1} h\underset{(3.13)}\to< c\kappa^{1/12}\kappa_1^{-1} \ve^{-1} |H|.
$$
Assume
$$
\ve_1< 10^{-3} \ve^3\tag 3.16
$$
and
$$
\kappa\ll \kappa_1^{24} \ve^{12}\tag 3.17
$$
($\leftrightarrow$ 3.15).

It follows that $|H_\beta \cap B|< \kappa_1 |H| $ for $\beta\not\in I_0$ while
certainly
$$
|H_{\alpha+\vp(y)} \cap B|> \frac 12 |H|.
$$
Hence $\alpha+\vp(y)\in I_0$ and we proved that
$$
I_0+\vp(y)=I_0 \text { in $G/H$ for all $y\in A\backslash A_1'$}.
$$
Thus
$$
\vp(A\backslash A_1')\subset \text { Sym}_1 (I_0) \text { in } G/H.\tag 3.18
$$
We distinguish two cases.

If $I_0=G/H$, then $|B|=|s(A_1')|> (1-\kappa_1)h'$.
We may then construct $A_1''$ as in \S2 and conclude (1.8) with $k< (2.4)$, $|\Omega|< \sqrt{\kappa_1} h'$.
\bigskip

Assume next $I_0\not= G/H$.
Hence Sym$_1(I_0)\not= G/H$ and $H'=\vp^{-1} \big(\text{Sym}_1 (I_0)\big) \supset H$ is a proper subgroup of $G$.
Hence
$$
\frac\ve 2 h'< |H'|\leq\frac{h'}2.
$$
By (3.18),
$$
A\backslash A_1'\subset H'.
$$
Since $I_0$ is a union of cosets of Sym$_1(I_0)$ in $G/H$, $\vp^{-1}(I_0)$ is a union of cosets $H_\tau'$ of $H'$, each satisfying
$$
|B\cap H_\tau'|>(1-\kappa_1) |H'| \text { for } \tau\in I_0'
$$ 
(by definition of $I_0$), where $I_0=\bigcup_{\tau\in I_0'} \text { Sym}_1(I_0)_\tau$.

Thus we may identify $H$ and $H'$ and write
$$
A\backslash A_1'\subset H
$$
with
$$
\frac \ve 2 h' <|H| <\frac {h'}2.\tag 3.19
$$
The set $s(A_1')=B_0\cup B_1$ with
$$
B_0=\bigcup_{\alpha\in I_0} \big(s(A_1')\cap H_\alpha) \ 
\text { and } \  B_1 =\bigcup_{\alpha\in I_1} \big(s(A_1')\cap H_\alpha)
$$
and
$$
\spreadlines{6pt}
\align
&|s(A_1')\cap H_\alpha|>(1-\kappa_1) |H| \ \text { for } \ \alpha \in I_0\tag 3.20\\
&|B_1| < c\kappa^{1/24} h'\tag 3.21\\
&I_0\not= \emptyset, I_0 \not= G/H.\tag 3.22
\endalign
$$
Next, take a set $z_1, \ldots, z_r \in A_1', r<\frac 2\ve$ of representatives for $\vp(A_1')$ and denote $A_2 =A_1'\backslash \{
z_1, \ldots, z_r\}$.
Then
$$
s(A_2)\subset s(A_1') \text { and } |s(A_2)|\geq 2^{-r} |s(A_1')|.
$$
Thus there is some $\alpha\in G/H$ s.t.
$$
|s(A_2)\cap H_\alpha|> \frac\ve 2 |s(A_2)|> \ve 2 ^{-r-1}|s(A_1')| >\ve 2^{-r-2} h'.
$$
Hence, for each $z\in s(z_1, \ldots, z_r)$
$$
|s(A_1')\cap H_{\alpha+\vp(z)}|\geq |\big(s(A_2) +z\big)\cap H_{\alpha+\vp(z)}| > \ve2^{-n-2}h'.
$$
We claim that $\alpha+\vp(z) =\beta\in I_0$. Otherwise, $\beta\in I_1$ and $s(A_1')\cap H_\beta\subset B_1$,
 implying by (3.21) that
$$
|s(A_1')\cap H_\beta|< c\kappa^{1/24 } h'
$$
and this is impossible, provided
$$
\kappa < 2^{-\frac{100}\ve}.\tag 3.23
$$
Hence
$$
I_0\supset \alpha+\vp \big(s(z_1, \ldots z_r)\big) = \alpha+\vp\big(s(A_1')\big)
$$ 
and since $I_0\subset\vp \big(s(A_1')\big)$, by (3.20), it follows that $I_0=\vp\big(s(A_1')\big)$ and therefore by (3.22)
$$
\vp\big(s(A_1')\big)\not= G/H.\tag 3.24
$$
Next partition
$$
I_0 =\vp\big(s(A_1')\big)= J\cup J'
$$
with
$$
J=\Big\{\alpha\in G/H, |A_1' \cap H_\alpha|>\frac{10}\ve\Big\}.
$$
Thus
$$
\Big|\bigcup_{\alpha\in J'} (A_1' \cap H_\alpha)|<\frac{20}{\ve^2}.\tag 3.25
$$
Take elements $\Cal Z =\big\{z_{\alpha, t}; \alpha\in J, t\leq \frac {10}\ve\big\} \cup\{z_\alpha; \alpha\in J'\}$ with
$\vp(z_{\alpha, t})=\alpha$. 

Then
$$
s(A_1')\supset s(\Cal Z)
$$
and
$$
\vp\big(s(A_1')\big)\supset \Big\{\sum_{\alpha\in J} u_\alpha\alpha; 0 \leq u_\alpha \leq \frac{10}\ve\Big\}+J'
=\langle J\rangle +J'
$$
where $\langle J\rangle$ is the group generated by $J\subset G/H$.
Thus $|\langle J\rangle| \leq |\vp\big(s(A_1')\big)|$.

>From (3.24), 
$ \langle J\rangle \not= G/H\text { and } H' =\vp^{-1} (\langle J\rangle)
$ is a proper subgroup of $G$.

Hence, by (3.25)
$$
|A_1'\backslash H'|<c(\ve)\tag 3.26
$$
and since $A\backslash A_1'\subset H$,
$$
|A\backslash H'|<c(\ve)
$$
with $H'$ a proper subgroup of $G$, $[G:H']\leq \frac 2\ve$.

Recalling the constraints (3.12), (3.15), (3.16), (3.17), (3.23) on the parameters $\ve, \ve_1, \kappa, \kappa_1$, take
$$
\spreadlines{6pt}
\align
\kappa_1& =\ve^2\\
\kappa&= 2^{-\frac{100}{\ve}}\\
\ve_1&= \Big(\frac 1\ve\Big)^{C\cdot 2^{\frac{100}{\ve}}}.
\endalign
$$

\bigskip
\noindent
{\bf (4).} Summarizing the preceding, we proved the following.

\proclaim
{Theorem 1}

Let $G$ be a finite abelian group and $A\subset G, |G|=h'$.
Let $\ve>0$ be a small constant.

There are the following alternatives.

\noindent
{\rm (4.1) } \ $s(A)=G$

\noindent
{\rm (4.2) } \ There is a proper subgroup $H$ of $G$, such that
$$
[G:H]< \frac 2\ve \text { and } \ |A\backslash H|< c(\ve).
$$
{\rm (4.3) } \ There are $k$ elements $x_1, \ldots, x_k\in A$ and a subset $\Omega_{x_1, \ldots, x_k} \subset G$ depending only on $x_1,
\ldots, x_k$, such that
$$
k<(1+\ve)\frac {\log h'}{\log 2} + c\log\log h'+ c(\ve)\tag 4.4
$$
$$
|\Omega_{x_1, \ldots, x_k}|\leq \ve h'+k\tag 4.5
$$
and
$$
A\subset \Omega_{x_1, \ldots, x_k}.\tag 4.6
$$
\endproclaim

\bigskip

\noindent  {\bf \S2.   Application to the class group}

\noindent
{\bf (5).} We apply the preceding to the class group $\Cal C$ for a large discriminant $D<0$.
\bigskip

Let $n\in\Bbb Z_+$ be square free; $n = \prod p_j$ with $(p_j, D)=1$ and $\Cal X_D(p_j)\not = -1$.
Let $C_j, C_j^{-1}$ be the classes that represent $p_j$.
Then $n$ is representable by all classes in the formal expansion $\prod\{C_j, C_j^{-1}\}$
(see [Bl], Cor. 2.3).

Let $G= \Cal C^2$.  Thus $h' =|G|=h/g$ with $g=|\Cal C/\Cal C^2|$ the number of genera.  Let $A= \{C_j^2\}\subset G$. We have
$$
\prod\{ C_j; C_j^{-1}\} =\Big(\prod C_j^{-1}\Big) s(A)\tag 5.1
$$
with $s(A)$ defined as in (1.0).

Fix $\ve>0$ a small parameter and apply Theorem 1 to $A\subset G$.

If $s(A)=G$ as in (4.1) of Theorem 1, then
$$
\prod \{C_j, C_j^{-1}\} =\Big(\prod C_j^{-1}\Big) \Cal C^2.
$$
Since $\Cal C/\Cal C^2$ is the group $\Cal G $ of  the genera, it follows
that in this case $n$ is representable by any form of the genus if it's representable by some form.

Assume now that $A$ satisfies the conditions of alternative (4.2) of Theorem 1.

Denote $\eta:\Cal C\to\Cal C^2$ obtained by squaring and let $\Cal C'= \eta^{-1} (H)$.  Since $\Cal C$ is a proper subgroup of $\Cal C_1$, we have
$$
\frac\ve 2 h<|\Cal C'|\leq\frac h2
$$
where $h=|\Cal C|$ is the class number.

There is a set of indices $\Cal J$ such that $|\Cal J|< C(\ve)$ and for $j\not\in\Cal J$, $C_j^2\in H$, hence $C_j, C_j^{-1}\in\Cal
C'$.

Denote $\Cal P_C$ the primes represented by the class $C$.
Thus $\Cal P_C =\Cal P_{C^{-1}}$.

It follows from the preceding that  $n \big(\prod_{j\in\Cal J} p_j\big)^{-1}$  
has all its prime factors in the set
$$
\Cal P(\Cal C')\equiv \bigcup_{C\in\Cal C'} \Cal P_c.
$$
We recall the following distributional theorem.

\proclaim
{Lemma 2} {\rm (Landau; [Bl], Lemma 5.1).}

Assume $D<(\log \xi)^A$, $A$ fixed.

Then
$$
|\{p\in\Cal P_C; p\leq\xi\}|=\pi_C(\xi) =\frac 1{\ve(C)h} \int^\xi_1\frac{dt}{\log t} +
\text{O}(\xi e^{-c\sqrt{\log \xi}})\tag 5.2
$$
with $\ve(C)=2$ if $C$ is ambiguous and $\ve(C)= 1$ otherwise.
\endproclaim

Recall also that the number of ambiguous classes equals
$$
\gamma = \# (\Cal C/\Cal C^2) \text { = number of genera $\ll 2^{\omega(D)}$}.
$$

Hence from (5.2)
$$
\spreadlines{6pt}
\align
\pi_{\Cal C'}(\xi)&= |\{p\in\Cal P(\Cal C'); p\leq\xi\}|\\
&\leq \sum_{C \text { ambiguous}} \  \pi_{C}(\xi) +\frac 12\sum_{\Sb C\in \Cal C'\\ \text { not ambiguous}\endSb}
\pi_{C}(\xi)\\
&\leq (\gamma +|\Cal C'|) \frac 1{2h} \int_2^{\xi} \frac {dt}{\log t} + \text{O}(\xi e^{-c\sqrt{\log \xi} } h)
\endalign
$$
and since $|\Cal C'| \leq \frac h2 \text { and } h<D^{\frac 12+\ve}< (\log \xi)^A$
$$
<\Big(\frac 14+ \frac 1{h^{1-\ve}}\Big) \int_2^{\xi} \frac {dt}{\log t}.\tag 5.3
$$

Thus, in summary, the number of integers $n\leq X$ obtained in alternative (4.2), is at most
$$
\sum_{\Sb  r\leq C_\ve; p_1\ldots p_r <X\\ \Cal C'<\Cal C\\ 2\leq [\Cal C; \Cal C']\leq \frac 2\ve\endSb}
\# \Big\{ n\leq \frac X{p_1\ldots p_r}; n \text { square free with prime factors in $\Cal P(\Cal C')$}\Big\}\tag 5.4
$$
with $\Cal P(\Cal C')$ satisfying (5.3) and $\{p_1, \ldots, p_r\}$ unordered and distinct, with $\Cal X_D(p_j) \not= -1$.

To bound the expressions $\# \{ \cdots \}$, use the upper bound sieve.

For instance Corollary 6.2 in [I-K], which we apply with $\Cal A=\Bbb Z_+$ and considering
$$
P(z)=\prod_{\Sb p\not\in\Cal P(\Cal C')\\ p< z\endSb} p.\tag 5.5
$$
Hence $g(d) =\frac 1d, |r_d (\Cal A)|\leq 1, \kappa =1, K=1$,
$$
V(z) =\prod_{p|p(z)} \big(1-g(p)\big) =\prod_{\Sb p<z\\ p\not\in \Cal P(\Cal C')\endSb} \Big(1-\frac 1p\Big)\tag 5.6
$$
and from [IK], (6.2), (6.80), applied with $D=z, s=1$
$$
\# \{ n<X; \big(n, p(z)\big) = 1\} < CXV(z)+R(z)\tag 5.7
$$
with
$$
R(z) =\sum_{d|P(z), d<z} |r_d(\Cal A)|\leq z.\tag 5.8
$$

Using (5.3) and partial summation
$$
\spreadlines{8pt}
\align
&\sum_{\Sb p<z\\ p\not\in\Cal P(\Cal C')\endSb}\frac 1p>
 \sum_{u< z} \frac 1{u^2} |\{ p\leq u; p\not\in \Cal P(\Cal C')\}|=
 \sum_{u<z} \frac 1{u^2} \Big(\frac u{\log u}-\pi_{\Cal C'} (u)\Big)+ \text{O}(1)\\
&> \sum_{\text{exp} (h^{1/A}) < u< z}\Big(\frac 34 -\frac 1{h^{1-\ve}}\Big) \frac 1{u\log u}+\text{O}(1)\\
&>\Big(\frac 34 -\frac 1{h^{1-\ve}}\Big) \log\Big(\frac{\log z}{h^{1/A}}\Big).\tag 5.9
\endalign
$$
Hence
$$
V(z) \lesssim \exp \left(-\Big(\sum_{\Sb p< z\\ p\not\in \Cal P(\Cal C')\endSb} \frac 1p\Big)\right) < \Big(\frac{h^{1/A}}{\log z}\Big)^{\frac 34
-o(1)}\tag 5.10
$$
for $z>\exp (D^{1/A})$.

Substituting in (5.7) with $z=\sqrt Y$ gives for $Y>\exp D^{1/A}$
$$
\# \{n<Y, n \text{ squarefree with prime factors in $\Cal P(\Cal C')$}\}
\lesssim \frac {h^{1/A}}{(\log Y)^{3/4-o(1)}} Y\tag 5.11
$$
(here $A$ is an arbitrary large fixed constant).

Returning to (5.4), we have for $\tau>0$ fixed, $X^\tau>\exp (D^{1/A})$
$$
\spreadlines{6pt}
\align
&\sum_{\Sb p_1\ldots p_r < X^{1-\tau}\\ \Cal X_D(p_j) \not= 1\endSb} 
\# \Big\{ n\leq\frac X{p_1\ldots p_r}; n \text { square free with primes in $\Cal P(\Cal C')$}\Big\}\\
&\overset{(5.11)}\to\lesssim \frac {h^{1/A} X}{\tau(\log X)^{3/4-o(1)}}\sum_{\Sb p_1\ldots p_r<X\\ \Cal X_D(p_j)\not= -1\endSb}
\frac 1{p_1\ldots p_r}\\
& \lesssim \frac {h^{1/A}X}{\tau(\log X)^{3/4 -o(1)}} \frac {(\frac 12\log\log X)^r}{r!}\\
&\lesssim \frac {h^{1/A} X}{\tau(\log X)^{\frac 34 -o(1)}}\tag 5.12
\endalign
$$
since $r< C(\ve)$.

This gives the contribution
$$
\ll\#\Big\{\Cal C'< \Cal C; [\Cal C:\Cal C']\leq \frac 2\ve\Big\}. \frac {h^{o(1)}X}{ \tau(\log X)^{\frac 34 -o(1)}}.
\tag 5.13
$$
It remains to consider the case (*): $n< X$ with prime divisors $p_1, \ldots, p_r$ such that $p_1\ldots p_r> X^{1-\tau}$.

\proclaim
{Lemma 3} Fix $r\in\Bbb Z_+$.
Then, for $X$ large enough
$$
|\{n<X; n \text { represented by $\Cal C$ and product of at most} \text { $r$ distinct primes}\}|
$$
$$
< \frac {rX}{\log X} \Big(\frac {e(\frac 12+\ve)\log\log X}{r-1}\Big)^{r-1}.\tag 5.14
$$
\endproclaim

\noindent
{\it Proof.}

We get the estimate
$$
\sum_{\Sb p_1< \cdots< p_{r-1}\\ \Cal X_D(p_j) \not= -1\endSb} \ \frac{X(p_1\cdots p_{r-1})^{-1}}
{\log(X(p_1\cdots p_{r-1})^{-1})}
$$
and since $ p_1\ldots p_{r-1} < X^{\frac {r-1}r}$, this is
$$
\spreadlines{6pt}
\align
&< \frac{rX}{\log X} \sum_{\Sb p_1< \cdots< p_{r-1} < X\\ \Cal X_D(p_j)\not= -1\endSb}
\ \frac 1{p_1 \cdots p_{r-1}}<\\
&\frac{rX}{(r-1)! \log X} \Big(\sum\limits_{\Sb p<X\\ \Cal X_D (p)\not= -1\endSb}\frac 1p \Big)^{r-1}.\tag 5.15
\endalign
$$
>From Lemma 2 and partial summation
$$
\spreadlines{8pt}
\align
&\sum_{\Sb p<X\\ \Cal X_D(p)\not= -1\endSb} \frac 1p =\frac 12 \sum_{C \text { non-ambiguous}} 
\ \sum_{\Sb p\in \Cal P_C\\ p< X\endSb} \frac 1p +\sum_{C \text { ambiguous}} \sum_{\Sb p\in\Cal P_C\\ p< X\endSb}\frac 1p\\
&\leq (h-|\Cal G|)\Big[\frac 1{2h}\int^X \frac {1}{y^2}\Big(\int^y_2 \frac {dt}{\log t}\Big) dy+c_A\int^X _{\exp (D^{1/A})}
\frac 1y \ e^{-c\sqrt{\log y}} dy\Big]\\
&+|\Cal G| \Big[\frac 1{2h}\int^X \frac 1{y^2} \Big(\int^y_2 \frac {dt}{\log t}\Big) dy + C_A\int^X_{\exp (D^{1/A})}\frac 1y
\ e^{-\sqrt{\log y}} dy\Big]\\
&+\int^{\exp (D^{1/A})}\frac 1{y^2} \ \frac {y}{\log y} dy\\
&< \frac 12 \int^X \frac 1{y^2} \Big(\frac y{\log y}+\text{O} \Big( \frac y{(\log y)^2}\Big)\Big) dy +C_A h \ e^{-c |D|^{(A/2)}}
+\frac {\log |D|}A\\
&<\frac 12 \log\log X+ \frac{\log |D|}A +\text{O}(1)\\
&<\Big(\frac 12+ o(1)\Big)\log\log X\tag 5.16
\endalign
$$
for $X$ large enough.

Substitution of (5.16) in (5.15) gives by Stirling
$$
\frac {rX}{\log X} \Big( \frac {e(\frac 12 +\ve)\log\log X}{r-1}\Big)^{r-1}
$$
proving Lemma 3.

Returning to the case $(*)$, we obtain the bound
$$
\spreadlines{8pt}
\align
&\sum_{\Sb u<X^\tau \\ \text {sq-free represented by $\Cal C$} \endSb} \ \frac {rXu^{-1}}{\log X} \Big(\frac {e(\frac 12+\ve) \log\log X}
{r-1}\Big)^{r-1}\\
&\sim \frac{rX}{\log X} \Big(\frac {e(\frac 12+\ve)\log \log X}{r-1}\Big)^{r-1} \Big\{\sum_{\Sb u< X^\tau\\ sq \text 
{ represented by $\Cal C$}\endSb} \frac 1{u}\Big\}\\
&\lesssim \frac {rX}{\log X} \Big(\frac { e(\frac 12+\ve)\log\log X}{r-1}\Big)^{r-1}
\Big(\int^{X^\tau}\frac {h^{1/A}}{u(\log u)^{\frac 12-o(1)}} +D^{1/A}\Big)\\
&\lesssim rX\Big(\frac{e(\frac 12+\ve) \log\log X}{r-1}\Big)^{r-1} \Big(\frac{\sqrt \tau}{(\log X)^{\frac 12-\frac 1A}}+
\frac 1{(\log X)^{1-\frac 1A}}\Big).\tag 5.17
\endalign
$$
Recall that $r<C_\ve$.
Taking $\tau=(\log X)^{-\frac 18}$, we obtain
$$
(5.13) +(5.17)< \# \Big\{\Cal C'<\Cal C; [\Cal C:\Cal C']\leq \frac 2\ve\Big\} \frac X{(\log X)^{\frac 12+\frac 1{32}}}
\tag 5.18
$$
with $\ve$ fixed, $D<(\log X)^C$ and $X$ large enough.
\bigskip

Next, consider the contribution from alternative (4.3) of Theorem 1.

This contribution is clearly bounded by
$$
\sum_{k \text{ as in }(4.4)} \ \sum_{\Sb p_1< \cdots < p_k\\ \Cal X_D(p_j)\not= -1\\ p_1\ldots p_k < X\endSb} \Big|\Big\{
n<\frac X{p_1\ldots p_k}; n \text { square free with primes in $\Cal P\big(\eta^{-1}(\Omega_{p_1, \ldots, p_k})\big)$
}\Big\}\Big|\tag 5.19
$$
where $\Omega_{p_1, \ldots, p_k} \subset\Cal C^2$ satisfies by (4.5)
$$
|\Omega_{p_1\ldots p_k}|< 2\ve |\Cal C^2|
$$
and hence $\tilde\Omega_{p_1\ldots p_k}=\eta^{-1} (\Omega_{p_1, \ldots, p_k})$ satisfies
$$
|\tilde\Omega_{p_1\ldots p_k}|< 2\ve h.\tag 5.20
$$
Repeating (5.3)-(5.11) with $\Cal C'$ replaced by $\tilde \Omega_{p_1, \ldots, p_k}$, we obtain that for \hfill\break
$y>\exp (D^{1/A})$
$$
\spreadlines{6pt}
\align
|\{ n<Y,  &\text{ $n$ square free with prime factors in $\Cal P(\tilde\Omega_{p_1\ldots p_k})\}|$}\\
& \lesssim \frac {h^{1/A}}{(\log Y)^{1-3\ve}} Y.\tag 5.21
\endalign
$$
We will consider several cases.

Assume in (5.19), $p_1\ldots p_k< \sqrt X$.

By (5.21), we obtain the bound
$$
\spreadlines{6pt}
\align
&\frac {h^{1/A}X}{(\log X)^{1-3\ve}} \ 
\sum_{\Sb p_1<\cdots< p_k< X\\ \Cal X_D(p_j)=-1\endSb} \ \frac 1{p_1\ldots p_k}\\
&<\frac {h^{1/A} X}{(\log X)^{1-3\ve}} \ \Big(\frac{\frac e2\log\log X}{k}\Big)^k.\tag 5.22
\endalign
$$
By (4.4), $k<(1+2\ve)\frac {\log h}{\log 2}$.
At this point, the size of $h$ becomes essential.
Write $h=(\log X)^\rho$ and $k=\sigma \log\log X$ with $\sigma <\frac{(1+2\ve)}{\log 2}\rho$.

Then (5.22) becomes
$$
\frac {h^{1/A} X}{(\log X)^{1-3\ve}} \Big(\frac e{2\sigma}\Big)^{\sigma\log\log X}
=\frac {h^{1/A} X}{(\log X)^{1-3\ve}} (\log X)^{1-\log 2\sigma)\sigma}.\tag 5.23
$$
Assume $\kappa>4\ve$ and 
$$
\rho< (1-\kappa)\frac{\log 2}2.\tag 5.24
$$
Then
$$
\sigma<\Big( 1-\frac \kappa 2\Big)\frac 12
$$
and hence
$$
\spreadlines{6pt}
\align
(5.23) &<\frac {h^{1/A} X}{(\log X)^{1-3\ve}} \ (\log X)^{\frac 12 (1-\log(1-\frac \kappa 2))(1-\frac \kappa 2)}\\
&<\frac X{(\log X)^{1-3\ve-\frac 1A}} \ (\log X)^{\frac 12 -\frac 1{16} \kappa^2+0(\kappa^3)}\\
&\leq \frac X{(\log X)^{\frac 12+\frac {\kappa^2}{20}}}\tag 5.25
\endalign
$$
provided
$$
\kappa> 10\Big(\sqrt \ve+\frac 1{\sqrt A}\Big).\tag 5.26
$$
Next, assume in (5.19) that $p_1\ldots p_k>\sqrt X$. Hence
$$
p_k>X^{\frac 1{2k}}.
$$
Rewrite the $p_1, \ldots, p_k$ sum in (5.19) as
$$
\sum_{\Sb p_1<\cdots< p_{k-1}\\ \Cal X_D(p_j)\not= -1\\ p_1\ldots p_{k-1} < X^{1-\frac 1{2k}}\endSb} \
\sum_{X^{\frac 1{2k}}<p_k<\frac X{p_1\ldots p_{k-1}}} \
\Big|\Big\{ n<\frac X{p_1\ldots p_k}; n \text { sf with primes in $\Cal P(\tilde\Omega_{p_1\ldots p_k})$}\Big\}\Big|.\tag 5.27 
$$
Fix $p_1\ldots p_k$, and denote $X'= \frac X{p_1\ldots p_{k-1}}> X^{\frac 1{2k}}$.

If $\frac X{p_1\ldots p_k} \leq \exp |D|^{1/A}$, then $\frac {X'}{\exp D^{1/A}}<p_k< X'$ and we obtain the bound 
$$
\spreadlines{6pt}
\align
X'\sum_{\frac {X'}{\exp (D^{1/A})}< p_k< X'} \frac 1{p_k} & < X' \{\log\log X' -\log (\log X' -D^{1/A})\}\\
& \lesssim X' \frac{D^{1/A}}{\log X'} \leq \frac {kX}{p_1\ldots p_{k-1}} \ \frac {|D|^{1/A}}{\log X}.\tag 5.28
\endalign
$$
If $\frac X{p_1\ldots p_k} > \exp |D|^{1/A}$, apply (5.21) to get the bound
$$
h^{1/A} X' \Big(\sum_{X^{\frac 1{2k}}<p_k<X'} \ \frac 1{(\log \frac{X'}{p_k})^{1-3\ve}} \ \frac 1{p_k}\Big).\tag 5.29
$$
If $2^\ell <\log \frac {X'}{p_k}< 2^{\ell+1}$, then $p_k>X'. e^{-2^\ell}$ and we obtain the contribution
$$
\spreadlines{6pt}
\align
&\lesssim \frac 1{2^{\ell(1-3\ve)} }\, [ \log\log X'-\log\log (X' e^{-2^\ell})]\\
&=-\frac 1{2^{\ell(1-3\ve)}} \, \log \Big(1-\frac{2^\ell}{\log X'}\Big).\tag 5.30
\endalign
$$
Distinguishing the cases $2^\ell <\frac 12\log X'$ and $\frac 12\log X'\leq 2^\ell <\log\frac {X'}{p_k}$ we get
$$
(5.30) <\frac {\log\log X}{(\log X')^{1-3\ve}} \lesssim \frac {k\log\log X}{(\log X)^{1-3\ve}}\tag 5.31
$$
and
$$
(5.29) < \frac {k|D|^{1/A} (\log\log X)}{(\log X)^{1-3\ve}} \ \frac X{p_1\ldots p_{k-1}}\tag 5.32
$$
that also captures (5.28).

Substitution of (5.32) in (5.27) gives the bound
$$
\spreadlines{6pt}
\align
&\frac {k|D|^{1/A}(\log\log X)}{(\log X)^{1-3\ve}} X \Big(\sum_{\Sb p_1<\cdots< p_{k-1}<X\\ \Cal X_D(p_j)=-1\endSb}
\ \frac 1{p_1\ldots p_{k-1}}\Big) \\
&<\frac {(\log\log X)^2 |D|^{1/A}} {(\log X)^{1-3\ve}} X \Big(\frac{\frac e2 \log\log X}{k-1}\Big)^{k-1}\tag 5.33
\endalign
$$
for which the bound (5.25) on (5.22) holds, under the assumption
$$
h<(\log X)^{(1-\kappa)\frac {\log 2}2}\tag 5.34
$$
with $\kappa$ satisfying (5.20).

In view of the preceding, in particular estimates (5.18) and (5.25), and taking into account that $|D|^{\frac 12-\ve}\ll h\ll
|D|^{\frac 12+\ve}$ and the number of genera is bounded by $2^{\omega(D)}\ll |D|^\ve$, we conclude

\proclaim
{Theorem 2}
Let $\kappa>0$ be a fixed constant and $D<0$ a negative discriminant satisfying
$$
|D|<(\log X)^{(1-\kappa)\log 2}.\tag 5.35
$$
Let $\Cal C$ be the class group.
Then for $X$ large enough
$$
\spreadlines{6pt}
\align
& \# \{n\sim X; n \text { square free, representable by some form but not by all forms of the genus}\}\\
&\lesssim_\kappa \# \Big\{ \text{$\Cal C'$ subgroup of $\Cal C; [\Cal C:\Cal C']<\frac{10^3}{\kappa^2}$}\Big\}. 
 \frac X{(\log X)^{\frac 12 +\frac {1}{33}}} +\frac X{(\log X)^{\frac 12+\frac{\kappa^2}{20}}}.\tag 5.36
\endalign
$$
\endproclaim

\bigskip
\noindent
{\bf \S3. Representation of shifted primes}

\noindent
{\bf (6).} Next, we establish a version of Theorem 2 for shifted primes.  

More precisely we get a bound on
$$
\spreadlines{6pt}
\align
\#\{q\sim X &\text { prime; $q+a$ square-free and representable by some form but not}\\
& \text{ all of the forms of the genus}\} \tag6.1
\endalign
$$

We use a similar strategy, combining the combinatorial Theorem 1 with upper bound sieving.
In fact, only the following crude upper bound will be needed.

\proclaim
{Lemma 4} Let $Y\in\Bbb Z$ be a large integer and let for each prime $\ell <Y$ a subset $R_\ell\subset \Bbb Z/\ell\Bbb Z$
be given, $|R_\ell|\in \{0, 1, 2\}$ (for the applications below). Then
$$
\#\{n < Y; \pi_\ell(n) \not\in R_\ell \text{  for each $\ell$}\}< 
(\log\log Y)^3 \prod_\ell \Big(1-\frac {|R_\ell|}\ell \Big) Y+\frac Y{(\log Y)^{10}}.\tag 6.2
$$
\endproclaim

\noindent
{\it Proof.}

Denote $\Omega=\{n\in\Bbb Z_+; n< Y\}$.

 For $\ell$ prime, let
$$
\Omega_\ell =\{n\in \Omega; \pi_\ell (n)\in R_\ell\}
$$
and from $m$ squarefree, denote
$$
\Omega_m =\bigcap_{\ell|m} \Omega_\ell.
$$
We have to bound
$$
\Big|\bigcap_\ell (\Omega\backslash \Omega_\ell)|\leq |\bigcap_{\ell< Y_0}(\Omega\backslash \Omega_\ell)|\tag 6.3
$$
with $Y_0<Y$ to be specified.

>From the inclusion-exclusion principle
$$
(6.3)\leq Y-\sum_{\ell< Y_0}|\Omega_\ell|+\sum_{\ell_1< \ell_2< Y_0} |\Omega_{\ell_1\ell_2}|\ldots
+\sum_{\ell_1<\cdots <  \ell_r<Y_0}| \Omega_{\ell_1\ldots\ell_r}|\tag 6.4
$$
with $r\in\Bbb Z_+$ even (to specify).

Clearly
$$
|\Omega_m|=\Big(\prod_{\ell|m}\frac{|R_\ell|}\ell\Big) Y+ 0\Big(\prod_{\ell|m} |R_\ell|\Big).\tag 6.5
$$
>From (6.4), (6.5)
$$
\spreadlines{6pt}
\align
\frac {(6.3)}Y
& \leq 1-\sum_{\ell<Y_0} \frac {|R_\ell|}\ell +\cdots + \sum_{\ell_1<\cdots <\ell_r<Y_0}
\frac {|R_{\ell_1}|}{\ell_1} \cdots \frac {|R_{\ell_r}|}{\ell_r}\\
&+\frac 1Y \Big(\sum_{\ell<Y_0} |R_\ell|+\ldots+ \sum_{\ell_1< \cdots< \ell_r <Y_0} |R_{\ell_1}|\cdots |R_{\ell_r}|\Big)\\
&\leq \prod_{\ell<Y_0}\Big(1-\frac {|R_\ell|}\ell\Big)+\sum_{r_1> r} \frac 1{r_1!} \Big(\sum_{\ell<Y_0} \ \frac
{|R_\ell|}{\ell}\Big)^{r_1} +\frac {2^{r+1}}Y \pmatrix Y_0+r\\ r\endpmatrix\\
&<\exp \Big( 3\sum_{\Sb Y_0<\ell<Y\\ \ell \text { prime}\endSb}\frac 1\ell\Big).\prod _{\ell<Y} \Big(1-\frac{|R_\ell|}\ell\Big)+
\sum_{r_1>r} \Big(\frac {2e\log\log Y}{r_1}\Big)^{r_1} +(3Y_0)^r Y^{-1}.
\endalign
$$
Take $r= 10^2\log\log Y, Y_0=Y^{10^{-3}(\log\log Y)^{-1}}$ to obtain (6.2).

\bigskip
Returning to Theorem 1 and alternative (4.2), we have 
$$
X\sim n=q+a =p_1\ldots p_r m \qquad \text { (square free)}\tag 6.6
$$
where $m$ has its prime factors in $\Cal P(\Cal C')$.
Let $X' =\frac X{p_1\ldots p_r}$.

Thus if $\ell\not\in \Cal P(\Cal C')$, $\pi_\ell (m)\not= 0$.
Also, since $q$ is prime, we have for any prime $\ell<\frac {X'}4, \ell \not= p_1, \ldots p_r$
$$
\pi_\ell (m) \not= \pi_\ell (a)/ \pi_\ell (p_1\ldots p_r).
$$
Hence we define for $\ell \in\Cal P(\Cal C'), \ell <\frac {X'}4, \ell \not= p_1,\ldots, p_r$
$$
R_\ell =\{ \pi_\ell(a)/\pi_\ell(p_1\ldots p_r)\},
$$
and for $\ell \not\in \Cal P(\Cal C'), \ell <\frac {X'}4, \ell \not= p_1, \ldots, p_r$
$$
R_\ell=\{0, \pi_\ell(a)/\pi_\ell (p_1\ldots, p_r)\}
$$
and $R_\ell =\phi$ otherwise.

Hence, recalling (5.3) and partial summation
$$
\spreadlines{6pt}
\align
\sum\frac {|R_\ell|}\ell &=\sum_{\Sb\ell\in \Cal P(\Cal C')\\ \ell <\frac {X'}4\endSb}\frac 1\ell +\sum_{\Sb \ell\not\in \Cal P(\Cal C')\\
\ell <\frac {X'}4\endSb} \frac 2\ell +O(\log\log r)\\
&= 2\log\log X' -\frac 14 \log\log X' +o(\log\log X')\tag 6.7
\endalign
$$
(for $\log X'> (\log X)^{1/A}$).

Therefore, given $p_1\ldots p_r$, the number of possibilities for $m$ in (6.6) is at most
$$
\frac {X'}{(\log X')^{7/4-}}\tag 6.8
$$
using (6.2).

Assume $X'>X^\tau$.
We obtain the bound (cf. (5.4))
$$
\spreadlines{6pt}
\align
&\# \Big\{\Cal C'< \Cal C; [\Cal C:\Cal C']\leq\frac 2\ve\Big\} \frac X{\tau^{7/4}(\log X)^{7/4-}} \ 
\sum_{\Sb p_1<\cdots< p_r<X\\ \Cal X_D(p_j)\not=-1, r< C(\ve)\endSb} \ \frac 1{p_1\ldots p_r}\\
&<\# \Big\{\Cal C'<\Cal C; [\Cal C:\Cal C']\leq\frac 2\ve\Big\} \frac X{\tau^{7/4}(\log X)^{7/4-}}\
\frac{(\frac e2 \log\log X)^r}r\\
&\ll_\ve \# \Big\{\Cal C'<\Cal C; [\Cal C:\Cal C']\leq \frac 2\ve\Big\}.\frac X{\tau^{3/4}(\log X)^{7/4-}}.\tag 6.9
\endalign
$$
For $X'<X^\tau$, proceed as follows.
Since $p_1<\cdots< p_r$ satisfies \hfill\break
$p_1\ldots p_r>\sqrt X$, we have $p_r>X^{\frac 1{2r}}$.

Writing
$$
n=q+a =p_1\ldots p_{r-1}. p_r . m
$$
and denoting $X'' =\frac X{p_1\ldots p_{r-1} m}$, fix $p_1, \ldots, p_{r-1}, m$
and estimate the number of possible $p_r\sim X''$.
Thus, for primes $\ell <\frac 14X'', (\ell, p_1\ldots p_{r-1}m)= 1$,
$$
\pi_\ell(p_r)\not\in \{0, \pi_\ell(a)/\pi_\ell (p_1\ldots p_{r-1} m)\}
$$
and, by Lemma 2, their number is at most
$$
(\log\log X)^4 \frac {X''}{(\log X'')^2} <\frac {r^2X''}{(\log X)^{2-}}.\tag 6.10
$$
This gives the contribution
$$
\sum_{\Sb p_1\ldots p_{r-1} m<X\\ m<X^\tau, m \text { sf}\\ \Cal X_D(p_1), \ldots, \Cal X_D(p_{r-1})\not= 
-1\\ \Cal X_D(p) \not=-1 \text
{ for } p|m\endSb} 
\ \frac X{p_1\ldots p_{r-1} m} . \frac 1{(\log X)^{2-}}\tag 6.11
$$
$$
\spreadlines{6pt}
\align
&<\frac X{(\log X)^{2-}}\Big(\frac{\frac e2 \log\log X}{r-1}\Big)^{r-1} \Big(\sum_{\Sb m<X^\tau\\ m  \text { sf, representable by
$\Cal C$} \endSb} \frac 1m\Big)\\
&\overset{(5.11)}\to < \frac X{(\log X)^{2-}} \Big(\frac{\frac e2 \log\log X}{ r}\Big)^r  h^{1/A}\big(\log (X^{\tau})\big)^{\frac 12
+}\\
& <\tau^{1/2}\frac X{(\log X)^{\frac 32-}}.\tag 6.12
\endalign
$$
Summing (6.9), (6.12) and appropriate choice of $\tau$, gives the bound
$$
\ll_\ve \# \Big\{\Cal C'<\Cal C; [\Cal C:\Cal C']\leq \frac 2\ve\Big\}\frac X{(\log X)^{\frac 32+\frac 1{20}}}\tag 6.13
$$
for the (4.2) exceptions.

Next the contribution of the (4.3) alternative from Theorem 1.

We have the bound (5.19), with the additional specification that $n=q+a$ ($q$ prime) and recalling that $\Omega_{p_1\ldots p_k}$ only
depends on the classes $C_1, \ldots, C_k\in \Cal C$ determined by $p_1, \ldots, p_k$.

Write again
$$
X\sim n=q +a=p_1\ldots p_{k-1} p_k m \text { with }  p_1<\cdots <p_k.
$$
Assuming $p_1\ldots p_k<\sqrt X$, we fix $p_1, \ldots, p_k$ and observe that the number of possibilities for $m$ with primes in 
 $\Cal P(\tilde\Omega_{ p_1,\ldots, p_k})$, is at most
$$
\frac X{p_1\ldots p_{k-1} p_k} .\frac 1{(\log X)^{2-4\ve}}\tag 6.14
$$
(using again Lemma 4).

This gives the contribution
$$
\frac X{(\log X)^{2-4\ve}} \ \sum_{\Sb p_1<\ldots< p_k\\ \Cal X_D(p_j)\not= -1\endSb} \ \frac 1{p_1\ldots {p_k}}\tag 6.15
$$
with $k\leq (4.4)$.
Following (5.22), (5.26), we get a bound
$$
\frac X{(\log X)^{\frac 32+\frac {\kappa^2}{20}}}\tag 6.16
$$
provided (5.24), i.e.
$$
\frac{\log h}{\log\log X} < (1-\kappa)\frac {\log 2}2.\tag 6.17
$$
If $p_1\ldots p_k \geq \sqrt X$, then $p_k> X^{\frac 1{2k}}$.

Proceed as follows.

Fix $p_1, \ldots, p_{k-1}$.
Then specify the class $\{C, C^{-1}\}$ of $(p_k)$ so that $\tilde\Omega=\tilde\Omega_{p_1, \ldots, p_k}$ is specified.
Take $m$ with prime factors in $\Cal P(\tilde\Omega)$. Finally estimate the number of primes $p=p_k<\frac X{p_1\ldots p_{k-1} m}$
satisfying the condition
$$
p\text { represented by $C$}\tag 6.18
$$
$$
\pi_\ell (p) \not=  \pi_\ell(a)/ \pi_\ell (p_1\ldots p_{k-1} m) \text { if  } (\ell, p_1\ldots p_{k-1}m)= 1, \ell <\sqrt X.\tag 6.19
$$

\proclaim
{Lemma 5} Let $Y<X$. Then
$$
|\{p<Y, p\text { satisfies (6.18), (6.19)}\}|\ll\frac {(\log\log X)}{h^{1-\ve}} \ \frac Y{(\log Y)^2}.\tag 6.20
$$
\endproclaim

Thus we have the bound
$$
\sum_{\Sb p_1< \ldots < p_{k-1}\\ \Cal X_D(p_j)\not=- 1\endSb} \ \sum_{C\in\Cal C}  \ 
\sum_{\Sb m  \text { sf with primes in }\Cal P(\tilde\Omega)\\ m<\frac{X^{1-\frac 1{2k}}}
{p_1\cdots p_{k-1}}\endSb} \# \Big\{p\lesssim \frac X{p_1\ldots p_{k-1} m};
(6.18), (6,19)\Big\}\tag 6.21
$$
and applying Lemma 5.
$$
\spreadlines{6pt}
\align
\# \Big\{p\lesssim \frac X{p_1\ldots p_{k-1} m}; (6.18), (6.19)\Big\} &\ll \frac {k^2(\log\log X)^4X}
{(\log X)^2 \, {p_1\ldots p_{k-1}m \, h^{1-\ve}}}\\
&\overset {(4.4)}\to\lesssim \frac {(\log\log X)^6X}{(\log X)^6 \, {p_1\ldots p_{k-1}m \, h^{1-\ve}}}.\tag 6.22
\endalign
$$
Next, by (5.21),
$$
\sum_{\Sb m<\frac X{p_1\ldots p_{k-1}}\\ \text { with primes in $\Cal P(\tilde \Omega)$}\endSb}\frac 1m\lesssim h^{1/A} (\log X)^{3\ve}.\tag
6.23
$$

This gives the bound (after summation over $C\in \Cal C$)
$$
\spreadlines{6pt}
\align
&\frac {h^{\ve+1/A} X}{(\log X)^{2-4\ve}} \ \sum_{\Sb p_1<\ldots < p_{k-1}< X\\ \Cal X_D(p_j)\not= -1\endSb}
\Big(\frac 1{p_1\ldots p_{k-1}}\Big)< \\
&\frac {X}{(\log X)^{2-5\ve-\frac 1A}} \Big(\frac{\frac e2 \log\log X}{k-1}\Big)^{k-1}.\tag
6.24
\endalign
$$
Using (6.24) and $k<\frac {1+\ve}{\log 2} \log h$, the assumption (6.17) will again ensure (6.16).

Hence from (6.13), (6.16) and the preceding, we can conclude

\proclaim
{Theorem 3}
Let $\kappa>0$ be a fixed constant and $D<0, D$ not a perfect square s.t.
$$
|D|<(\log X)^{(1-\kappa)\log 2}.\tag 6.25
$$
Let $\Cal C$ be the class group.
Then, for $X$ large enough and $a\in\Bbb Z_+$, $a=o(X)$ fixed, we have

\noindent
$\#$ \{$q+a\sim X$; $q$ prime, such that $q+a$ is squarefree and representable by some form but not by all forms of the genus\} $\lesssim_\kappa$.
$$
\# \{\Cal C' \text { subgroup of $\Cal C$;} [\Cal C:\Cal C']<C(\kappa)\}\frac X{(\log X)^{\frac 32+\frac 1{20}}}
+ \frac X{(\log X)^{\frac 32+\frac {\kappa^2}{20}}}.\tag 6.26
$$
\endproclaim

Note that by decomposing $\Cal C$ into cyclic groups, one easily gets a bound
$$
\# \{\Cal C' \text { subgroup of } \Cal C; [\Cal C:\Cal C']< C(\kappa)\}<
C_1 (\kappa) (\log h)^{C_1(\kappa)} < (\log\log X)^{C(\kappa)}.
$$

\noindent
{\bf Proof of Lemma 5.}

In order to estimate the size of the set
$$
\{p<Y, p \text { satisfies (6.18), (6.19)}\}\tag 6.27
$$
we factor in prime ideals and consider the larger set
$$
\{\alpha \in I; \alpha\in C, N(\alpha) <Y \text { and } \pi_\ell \big(N(\alpha)\big)\not\in R_\ell \text { for } \ell <Y_0\}\tag 6.28
$$
where $I$ denotes the integral ideals in $O_K, K=\Bbb Q(\sqrt D), D=D_0 f^2$ with $D_0<0$ squarefree, $N(\alpha)$ stands for the norm
of $\alpha$ and $\ell$ runs over primes,
$$
(6.29) \  \left\{\aligned R_\ell&= \{0, \xi_\ell\}, \xi_\ell =\pi_\ell (a)/\pi_\ell (p_1\ldots p_{k-1} m) \text { if } (\ell, p_1\ldots p_{k-1} m)= 1\\
R_\ell &= \{0\} \text { otherwise.}
\endaligned
\right.
$$
In fact, we restrict ourselves in (6.28) to primes $\ell <Y$ such that
$$
(\ell, p_1\ldots p_{k-1} m)=1.\tag 6.30
$$ 
Define
$$
\Omega =\{\alpha\in I; \alpha\in C, N(\alpha)<Y\}
$$
and
$$
\Omega_\ell =\{\alpha\in\Omega; \pi_\ell \big(N(\alpha)\big) \in R_\ell\}
$$
for $\ell$ prime,
$$
\Omega_n=\bigcap_{\ell|n}\Omega_\ell
$$
for $n$ square-free.

Proceeding as in the proof of Lemma 4, estimate
$$
\align
&\Big|\bigcap_{\ell<Y_0, (6.30)} (\Omega\backslash\Omega_\ell)\Big| \leq\\
&|\Omega|-\sum_{\ell<Y_0, (6.30)} |\Omega_\ell|+ \sum_{\ell_1<\ell_2< Y_0} |\Omega_{\ell_1, \ell_2}| -\cdots +\sum_{\ell_1< \cdots< \ell_r< Y_0}
|\Omega_{\ell_1, \ldots, \ell_r}|\tag 6.31
\endalign
$$
with $r\in\Bbb Z_+$ $r\sim\log\log Y$ suitably chosen.

We evaluate $|\Omega_n|$ using Hecke characters.

The condition that $\alpha\in C$ becomes
$$
\frac 1h \sum_{\lambda\in \widehat{\Cal C}} \overline{\lambda(C)} \lambda(\alpha)=1\tag 6.32
$$
where $\lambda$ runs over the class group characters $\widehat{\Cal C}$.

Denote $\Cal X_\ell$ the principal character of $\Bbb Q(\mod \ell)$.

If (6.30), $\pi_\ell\big(N(\alpha)\big)\in R_\ell$ may be expressed as
$$
1-\Cal X_\ell \big(N(\alpha)\big) +\frac 1{\ell-1}\sum_{\Cal X(\mod \ell)} \overline{\Cal X} (\xi_\ell) \Cal X \big(N(\alpha)\big)=1.\tag 6.33
$$
Thus
$$
|\Omega_n|=\sum_{N(\alpha)<Y} \Big[\frac 1h\sum_{\lambda\in\hat{\Cal C}} \overline{\lambda(C)} \lambda(\alpha) \Big] \prod_{\ell|n} (6.33).\tag 6.34
$$

We will use the following classical extension of the Polya-Vinogradov inequality for finite order Hecke characters.

\proclaim
{Proposition 6}

(i) Let $\Cal X$ be a non-principal finite order Hecke character $(\mod f)$ of $K$. Then
$$
\Big|\sum_{N(\alpha)< x} \Cal X(\alpha)\Big| < C\big (|D|N(f)\big)^{1/3} [\log |D|N(f)]^2 \ x ^{1/3}\tag 6.35
$$
and also

(ii)
$$
\sum_{N(\alpha)<x} 1= c_1x +0\big(|D|^{1/3} (\log |D|)^2\big) x^{1/3}\tag 6.36
$$
where
$$
c_1 =\prod_{p|f} \Big(1-\frac 1p\Big) L(1, \Cal X_D\Big).\tag 6.37
$$
\endproclaim

This statement follows from [L], (1), (2) p. 479; for (6.37), see [Bl], (2.5).

Analyzing (6.33), (6.34) more carefully, we see that
$$
\spreadlines{6pt}
\align
|\Omega_n|&= \frac 1h \sum_{N(\alpha)\leq Y} \prod_{\ell|n} \Big( 1-\frac {\ell-2}{\ell-1}\Cal X_\ell \big(N(\alpha)\big)\Big)\tag 6.38\\
& +0(6.39)
\endalign
$$
where (6.39) is a bound on sums
$$
\sum_{N(\alpha)<Y}\Cal X(\alpha) \text { with } \Cal X(\alpha) =\lambda(\alpha)\Cal X'\big(N(\alpha)\big)\tag 6.40
$$
where $\lambda\in\hat{\Cal C}, \Cal X'$ is a $(\mod n_1)$-Dirichlet character with $n_1|n$ and either $\lambda$ or $\Cal X'$ non-principal.
Thus by (6.35)
$$
(6.39)< C|D|.n Y^{1/3}< C|D| Y_0^r Y^{1/3}\tag 6.41
$$
which collected contribution in (6.31) is at most
$$
C|D|Y_0^{2r} Y^{1/3} <Y^{1/2}\tag 6.42
$$
imposing the condition
$$
|D|Y_0^r < Y^{\frac 1{20}}.\tag 6.43
$$

Analyzing further (6.37) using (6.36), we obtain
$$
\spreadlines{8pt}
\align
|\Omega_n|&=\frac {c_1} h. Y \prod_{\Sb \ell|n\\ \Cal X_D(\ell)=1\endSb} \Big[ 1-\frac {\ell-2}{\ell-1}\Big(1-\frac 1\ell\Big)^2\Big]
.\prod_{\Sb\ell|n\\ \Cal X_D(\ell)=0\endSb} \Big[1-\frac{\ell-2}{\ell-1}\Big(1-\frac 1\ell\Big)\Big] .
\\
&\prod_{\Sb \ell |n\\ \Cal X_D(\ell)=-1\endSb} \Big[1-\frac {\ell-2}{\ell-1}\Big(1-\frac  1{\ell^2} \Big)\Big]. \tag 6.44\\
&+ 0(Y^{1/2}).
\endalign
$$
Substituting (6.44) in (6.31) gives
$$
\spreadlines{6pt}
\align
&\frac {c_1}h Y \prod_{\Sb \ell < Y_0, (6.30)\\ \Cal X_D(\ell)=1\endSb} \Big( 1-\frac 3\ell+\frac 2{\ell^2}\Big). \prod_{\Sb \ell< Y_0 (6.30)\\ \Cal
X_D(\ell)=0\endSb} \Big(1-\frac 2\ell\Big). \prod_{\Sb \ell <Y_0, (6.30)\\ \Cal X_D(\ell) =-1\endSb} \Big(1-\frac 1\ell-\frac 2{\ell^2}\Big)+\\
& 0\Big(Y\frac 1{r!} \Big(\sum_{\Sb \ell < Y_0\\ \ell \text { prime}\endSb} \frac 3\ell\Big)^r
 +Y^{1/2} Y_0^r\Big)\\
&\ll \frac{|D|^\ve} h Y\frac{(\log\log X)^3}{(\log Y_0)^2} +0\Big(Y\Big(\frac {3\log\log Y_0} r\Big)^r +Y^{1/2} Y_0^r\Big).\tag 6.45
\endalign
$$
Taking $r=10^2 \log\log Y, Y_0 = Y^{10^{-4} (\log\log Y)^{-1}}$, (6.43) holds and we obtain (6.20).
This proves Lemma 5. \qed

Theorem 3 may be combined with Iwaniec' result [I] on representing shifted primes by the genus of a binary quadratic form (see the Appendix
for a quantitative review of that argument, when the quadratic form $Ax^2+ Bxy+Cy^2=f(x, y)$ is not fixed).
Thus, fixing $a\not= 0$, and assuming $D=B^2-4AC$ not a perfect square, it follows from [I] that
$$
\spreadlines{6pt}
\align
&\# \{q+a\sim X; q\text { prime and $q+a$ squarefree and representable by the genus of $f$}\}\\
&\gg \frac X{(\log X)^{3/2 +\ve}}\tag 6.46
\endalign
$$
and this statement is certainly uniform assuming $|A|, |B|, |C|< \log X(?)$

\proclaim
{Corollary 4} Let $f$ be as above with discriminant $D<0$,  and assume for some $\kappa>0$
$$
|D|<(\log X)^{(1-\kappa)\log 2}\tag 6.47
$$
with $X$ sufficiently large.  Then
$$
\spreadlines{6pt}
\align
&\# \{ q+a\sim X; q \text { prime, such that $q+a$ is representable by $f$\}}\\
&\gg \frac X{(\log X)^{3/2 +\ve}}.
\endalign
$$
\endproclaim

\vfill\eject

\noindent
{\bf Appendix}
\medskip

Let $\phi(x,y)$ be a primitive positive definite binary quadratic form of discriminant $-D$ where $D<\log X$, and let
$$
S_1(X,\phi,a)=\sum_{\Sb{p\leq X, p\not\; | D}\\ {p=f(x, y)+a} \\ {{(x,y)=1,f\in R_{\phi}}}\endSb}1
$$
where $R_{\phi}$ denotes the genus of $\phi$.  Then Theorem 1 of [Iw] gives us the following lower bounds for $S_1$.

\medskip
\noindent
{\bf Theorem A.1.}
%\begin{thm}\label{iwgeneral}
{\sl For $a\in\Bbb Z$ and $\phi$ a primitive positive definite binary quadratic form of discriminant $-D$ where $D\leq\log X$, let $S_1(X,\phi,a)$ be as above.
Then for $\epsilon>0$ we have
$$
S_1(X,\phi,a)\gg_{\epsilon} \frac{X\cdot D^{-\epsilon}}{(\log X)^{3/2}}
$$
where the implied constant does not depend on $D$.
}

The following two lemmas are essentially Theorems 2 and 3 from [Iw] in the case $D<\log X$, where the integer $m$ represented by $R_{\phi}$ is assumed to be square free and $(m,D)\leq 2$.

\medskip
\noindent
{\bf Lemma A.2.}
(Iwaniec). {\sl Let $-D<0$ be the discriminant of $f(x,y)=Ax^2+2Bxy+Cy^2$, and write
$$-D= -2^{\theta_2}\cdot p_1^{\theta_{p_1}}\cdots p_r^{\theta_{p_r}}, \quad D_p=p^{-\theta_p}\cdot D,$$
where $\theta_p\geq 1$ for $1\leq i\leq r$, and $\theta_2\geq 0$.  Write $m=\delta n=2^{\epsilon_2}n$ where $m$ is a positive square free integer (so $0\leq\epsilon_2\leq 1$) such that $(n,2D)=1$.
Then $m$ is represented by the genus of $f$ iff the conditions on $m$ in Table 1 are satisfied
\footnote{Table 1 also specifies a quantity $\kappa$ and $\tau$ for each described case.  These do not have to do with whether $m$ is represented or not, but will be used later.}.
}

With the notation above, for $p\not=2$, let
$$
\Cal L_p'(n)=\left\{l\; |\; 0<l<p, \left(\frac{l}{p}\right)=\left(\frac{A\cdot 2^{\epsilon_2}}{p}\right)\right\},
$$
$$
\Cal L_p''(n)=\left\{l\; |\; 0<l<p, \left(\frac{l}{p}\right)=\left(\frac{-A\cdot 2^{\epsilon_2}\cdot k(-D_p)}{p}\right)\right\}
$$

\bigskip
\vfill\eject
\centerline
{TABLE 1. Representation of $2^{\varepsilon_2} n$ by $f$}

\medskip
\vbox{\offinterlineskip
\hrule
\halign{&\vrule#&
\strut \, \hfil# \ \hfil\cr
height5pt&\omit&&\omit&&\omit&&\omit&&\omit&&\omit&\cr
& Description of $\theta_p$ \hfil&& $\Cal K$ && $\tau$ && {Contributions on $n$}&& {Contributions on $D$}
&& \ &\cr
\noalign{\hrule}
height 5pt&\omit&&\omit&&\omit&&\omit&&\omit&&\omit&\cr
& $\theta_{p_i}\geq 1, \; p_i\not=2$ && $\frac{p_i-1}{2}$ && $p_i$ && $\left(\frac{n}{p_i}\right) =
\left(\frac{A\cdot2^{\epsilon_2}}{p_i}\right)$ && none && (1)&\cr
\noalign{\hrule}
height 5pt&\omit&&\omit&&\omit&&\omit&&\omit&&\omit&\cr
& $p|m$, \; $\theta_p=0$ && 1 && 1 && none && $\left(\frac{-D}{p}\right)=1$ && (2)&\cr
\noalign{\hrule}
height 5pt&\omit&&\omit&&\omit&&\omit&&\omit&&\omit&\cr
&$\epsilon_2=0$, $\theta_2=0$ && $1$ && $1$ && none && $D\equiv -1 \; (4)$ && (3)&\cr
\noalign{\hrule}
height 5pt&\omit&&\omit&&\omit&&\omit&&\omit&&\omit&\cr
&$\epsilon_2=0$, $\theta_2=2$ && $1$ or $2$ && $4$ && $n\equiv A \; (4)$ or $n\equiv -A\, D_2 \; (4)$ &&
$D_2\equiv -1 \; (4)$ or  $D_2\equiv 1 \; (4)$ && (4)&\cr
\noalign{\hrule}
height 5pt&\omit&&\omit&&\omit&&\omit&&\omit&&\omit&\cr
&$\epsilon_2=0$, $\theta_2=3$ && $2$ && $8$ && $n\equiv A \; (8)$ or $n\equiv A\,(1-2D_2) \; (8)$ && none && (5)&\cr
height 5pt&\omit&&\omit&&\omit&&\omit&&\omit&&\omit&\cr
\noalign{\hrule}
height 3pt&\omit&&\omit&&\omit&&\omit&&\omit&&\omit&\cr
& $\epsilon_2=0$,  $\theta_{2}=4$ && $1$ && $4$ && $n\equiv A \; (4)$ && none && (6)&\cr
\noalign{\hrule}
height 5pt&\omit&&\omit&&\omit&&\omit&&\omit&&\omit&\cr
&$\epsilon_2=0$,  $\theta_{2}\geq 5$ && $1$ && $8$ && $n\equiv A\; (8)$ && none && (7)&\cr
\noalign{\hrule}
height 5pt&\omit&&\omit&&\omit&&\omit&&\omit&&\omit&\cr
&$\epsilon_2=1$,  $\theta_2=0$ && $1$ && $1$&& none && $D\equiv -1 \; (8)$ && (8)&\cr
\noalign{\hrule}
height 5pt&\omit&&\omit&&\omit&&\omit&&\omit&&\omit&\cr
&$\epsilon_2=1$, $\theta_2=2$ && $1$ && $4$ && $n\equiv A\,\frac{1-D_2}{2} \; (4)$ && $D_2\equiv -1 \; (4)$&& (9)&\cr
height 5pt&\omit&&\omit&&\omit&&\omit&&\omit&&\omit&\cr
\noalign{\hrule}
height 5pt&\omit&&\omit&&\omit&&\omit&&\omit&&\omit&\cr
&$\epsilon_2=1$, $\theta_2=3$ && $2$ && $8$ && $n\equiv -A\, D_2 \; (8)$ or $n\equiv A\,(2-D_2) \; (8)$ && none &&
(10)&\cr}
\smallskip
\hrule}
\smallskip
\bigskip

\noindent
where $k(-D_p)$ denotes the square free kernel of $-D_p$.  Note that each of $\Cal L_p'$ and $\Cal
L_p''$ always contains $(p-1)/2$ elements.  Define $\Cal L_2(n)$ as follows:
$$
\Cal L_2(n) =
\left\{
\aligned
\{l\; |\; 0<l<4, \, l\equiv A \; (4)\text { or } l\equiv -A\,D_2\; (4)\} & \text{if } \epsilon_2=0,\, \theta_2 = 2
\\\{l\; |\; 0<l<8, \, l\equiv A \; (8) \text{ or } l\equiv A\,(1-2\,D_2)\; (8)\} & \text{if } \epsilon_2= 0, \; \theta_2=3
\\\{l\; |\; 0<l<4, \, l\equiv A \; (4)\} & \text{if } \epsilon_2=0,\, \theta_2=4
\\\{l\; |\; 0<l<8, \, l\equiv A \; (8)\} & \text{if } \epsilon_2=0, \theta_2\geq 5
\\\{l\; |\; 0<l<8, \, l\equiv -A\,D_2 \; (8) \text{ or } l\equiv A\,(2-D_2) \; (8)\} & \text{if } \epsilon_2=1, \theta_2=3
\\\{l\; |\; 0<l<4, \, l\equiv -A\,\frac{D_2-1}{2} \; (4)\} & \text{if } \epsilon_2=1, \theta_2=2, \;D_2\equiv -1 \; (4)
\\\{0\} & \text{if } \epsilon_2\geq\theta_2.\\
\endaligned\right.
$$
Note that $\Cal L_2(n)$ contains $\kappa$ elements, where $\kappa$ is as in Table~1.  With this notation, we have

\bigskip
\noindent
{\bf Lemma A.3.}
(Iwaniec). {\sl Let $D$, $\theta_p$, $m$, $n$, and $\delta$ be as in Lemma 0.2, and let $\tau_2$ be the corresponding value of $\tau$ in the case $p=2$
in Table~1.  Define $Q=\tau_2\cdot\prod_{p_i|D_2}p_i$, and let
$$
P=\left\{p\; \Big|\; \left(\frac{k(-D)}{p}\right)=1\right\}
\tag A.1
$$
where $k(-D)$ is the square free kernel of $-D$. 
Then $m=2^{\epsilon_2}n$ is represented by the genus of $\phi$ iff $m$ satisfies the conditions in Table 1, all the prime factors of $n$ belong to $P$, and
$$
n\equiv L\; (Q)
$$
where $L>0$ is an integer satisfying the conditions

\roster
\item "{$\bullet$}" $0<L<Q$,
\item "{$\bullet$}" $L\equiv l \; (\tau_2)$ for some $l\in\Cal L_2(n)$,
\item "{$\bullet$}" for each $p_i|D_2$ there exists $l\in\Cal L'_{p_i}(n)$ such that $L\equiv l \; (p_i)$.
\endroster
\bigskip

Furthermore, if $\Cal L$ denotes the set of $L$ satisfying these conditions,
$\left(\frac{k(-D)}{L}\right)=1$ for each $L\in\Cal L$.
}

\noindent Let $\Cal P = \{\text{primes } p\not|D \; \text{s.t.} \; \left(\frac{k(-D)}{p}\right)=-1\}$, let $E= Q\delta$, and let $\phi_E(N)=\phi(N\cdot E)/\phi(E)$.  For $D$ fixed, it is crucial to
the $\frac{1}{2}$-dimensional sieve that the condition
$$
\left|\sum_{\Sb{p\leq z}\\{p\in\Cal P}\endSb}\frac{\log p}{\phi_E(p)}-\frac{1}{2}\log z\right|<c
\tag A.2
$$
is satisfied for some constant $c$ for all $z>1$.  In our case of $D\leq\log X$, this holds in the following form for some constant $C_1$ not depending on $D$:
$$
\left|\sum_{\Sb{p\leq z}\\{p\in\Cal P}\endSb}\frac{\log p}{\phi_E(p)}-\frac{1}{2}\log z\right|\ll_{\epsilon} C_1 D^{\epsilon}
\tag A.3
$$
for any $z\geq 1$.
This can be seen from the proof of Theorem 3.2.1 of \cite{Gl} and the fact that
$$
\sum_{\left(\frac{k(-D)}{p}\right)=1, p\leq z}\frac{\log p}{p} = \frac{\log z}{2} + D^{\epsilon}\cdot O(1)
$$
where the implied constant depends only on $\epsilon$.  As in [Iw], let
$$
C_0:=\lim_{z\rightarrow\infty}\prod_{\Sb{p<z}\\ {p\in\Cal P}\endSb}\left(1-\frac{1}{\phi_E(p)}\right)\sqrt{\log z}
$$
for which Iwaniec shows in [Iw]

\bigskip
\noindent
{\bf Lemma A.4.}
(Iwaniec). {\sl  Let $C_0$ be as above.  We have
$$
C_0 = e^{-\gamma/2}\prod_{\Sb{p\not{ \;}  |\, a}\\{p\in \Cal P}\endSb}\left(1-\frac{1}{(p-1)^2}\right)\cdot \prod_{p|\, Da}\left(1-\frac{1}{p}\right)^{-1/2}\cdot
\prod_{p\not{ \;} |\, Da}\left(1-\frac{1}{p}\right)^{-\left(\frac{-k(D)}{p}\right)/2}
$$}
\bigskip

Finally, we recall the following theorem of Bombieri:
\bigskip

\noindent
{\bf Lemma A.5.}
(Bombieri). {\sl Let $\pi(x,k,l)$ denote the number of primes less than $x$ which are $l$ modulo $k$.
There exists an absolute constant $U$ such that
$$\sum_{k<\frac{\sqrt{x}}{(\ln x)^U}}\max_{\Sb{l}\\ {(l,k)=1}\endSb}\left|\pi(x,k,l)-\frac{\text{\rm Li}\, x}{\phi(k)}\right|\ll\frac{x}{(\log x)^{20}}.$$
}
\bigskip

We are now ready to introduce the notation relevant to our problem and recall the lemmas resulting from the $\frac{1}{2}$-dimensional sieve.  For $L$ and $\delta$ as above, and $1<s\leq \frac{4}{3}$,
\medskip
\roster
\item "{$\bullet$}" $D_1=2\text{ or } 1=\text{ greatest divisor of } 2D \text{ prime to } Q\cdot a$
\item  "{$\bullet$}" $M=\{m\in \Bbb N \; | \; m=\frac{p-a}{\delta}, p\leq X, p\equiv \delta L+a\; (Q\delta), \; (m,D_1)=1\}$
\item  "{$\bullet$}" $M_d=\{m\in M \; | \; m\equiv 0 \; (d)\}$
\item  "{$\bullet$}" $Y= \phi(E)\cdot |M| =\text{\rm Li}\, (X)$
\item "{$\bullet$}" $R_d(M) = |M_d|-\frac{Y}{\phi(dE)}$
\item "{$\bullet$}" $y =\frac{\sqrt{X}}{Q\delta D(\log X)^U}$
\item "{$\bullet$}" $A(M,y^{1/s}) = \#\{m\in M \, \text{s.t. } \, m\not\equiv 0 \;(p),\; y^{1/s}>p\in\Cal P\}$
\endroster
\medskip

By Lemma~A.3, the following is precisely what is needed to evaluate $S_1$:
$$
\aligned
\sum_{\Sb{|a|<f(x,y)+a=p\leq X}\\ {(x,y)=1,f\in R_{\phi}}\endSb}1&=\sum_{d}\sum_{L\in\Cal L}\sum_{\Sb{X\geq p\equiv \delta L +a \; (Q\delta)}\\{q|((p-a)/\delta)\Rightarrow q\in P}\\ {((p-a)/\delta,2D)=1, p>|a|}\endSb}1\\
&=\sum_{\Sb{\delta}\\ {2|a\delta}\endSb}\sum_{\Sb{L\in \Cal L}\\ {(\delta L+a,Q\delta)=1}\endSb}\sum_{\Sb{m\in M}\\ {q|m\Rightarrow q\in P}\endSb}1 + \Cal R
\endaligned
\tag A.5
$$
where $\Cal R\leq 2|D|$.  It is the innermost sum in (A.5) that we evaluate with the help of the sieve.  Note that if $(d, QA)=1$, there exists an integer $d'$ such that $d'Q+L\equiv 0 \; (d)$ and
$$M_d=\{m\; | \; m=\frac{p+a}{\delta}, \; p\leq X \; p\equiv A + \delta\, L + Q\delta\, d' \; (Q\delta\, d)\}.$$
>From [Iw] we then have
$$
\left| |M_d| -\frac{\text{\rm Li} X}{\phi(Q\,\delta d)}\right|\leq 2 \max_{\Sb{l}\\ {(l,Q\,\delta\, d)=1}\endSb
}\left|\pi(X,Q\delta\,d,l)-\frac{\text{\rm Li} X}{\phi(Q\,\delta\, d)}\right|
\tag A.6
$$

With the notation above, the expression in (A.3) combined with the $\frac{1}{2}$-dimensional sieve gives the following in our case:

\bigskip
\noindent
{\bf Theorem A.6.}
{\sl $$
\aligned
A(M,y^{1/s})&\gg_{\epsilon} \sqrt{\frac{e^{\gamma}}{\pi}}\cdot \frac{C_0Y}{\phi(E)\sqrt{\log
3y}}\cdot\left(\int_1^s\frac{dt}{\sqrt{t(t-1)}}-\frac{(\log X)^{\epsilon}}{(\log 3y)^{1/10}}\right)-
\sum_{\Sb{d<y}\\ {p|d\Rightarrow p\in\Cal P}\endSb}|R_d(M)|\\
&\geq \frac{C_0\cdot\sqrt{2\frac{e^{\gamma}}{\pi}}}{\phi(Q\delta)}\cdot\frac{X}{(\log X)^{3/2}}\cdot\left(\int_1^s\frac{dt}{\sqrt{t(t-1)}}+(\log X)^{\epsilon}\cdot
o(1)\right)+\text{\rm O}(X\log^{-20} X).
\endaligned
$$}

The estimation of the remainder term comes from Lemma~A.5 and (A.6).  Also, for sufficiently large $X$ (such that $(X^{1/2}(\log X)^{-15-U})^{1/s}>X^{1/3}$) and $1<s<\frac{4}{3}$ we have
$$
A(M,y^{1/s})=\sum_{\Sb {m\in M}\\ {q|m\Rightarrow q\in P}\endSb}1+
\sum_{\Sb{p_1p_2m\in M}\\ {q|m\Rightarrow q\in P}\\ {{y^{1/s}\leq p_1,p_2\in\Cal P}}\endSb}1.
$$

We have a lower bound for $A(M,y^{1/s})$, and we would like a lower bound for the first sum in the equation above.  To this end, Iwaniec shows:

\bigskip
\noindent
{\bf Lemma A.7.}
(Iwaniec).
{\sl Let $|Q\delta|\ll(\log X)^{15}$ and $s>1$.  Then
$$
\sum_{\Sb{p_1p_2m\in M}\\ {q|m\Rightarrow q\in P}\\ {{y^{1/s}\leq p_1,p_2\in\Cal P}}\endSb}1<{\frac{4e^{\gamma/2}C_0\sqrt{s-1}}{\sqrt{\pi}\phi(Q\delta)\sqrt{s}}}\log(2s-1)\frac{4s^2X}{(\log
X)^{3/2}}(1+o(1))$$
}
\bigskip

Together with Theorem~A.6, for $1<s<\frac{4}{3}$ and $Q\delta\ll\log X$, this gives us
$$
\sum_{\Sb{m\in M}\\{q|m\Rightarrow q\in P}\endSb}1\gg
\sqrt{\frac{2e^{\gamma}}{\pi}}\cdot\frac{C_0}{\phi(Q\delta)}\cdot\frac{X}{(\log
X)^{3/2}}\cdot\left(\int_1^s\frac{dt}{\sqrt{t(t-1)}}-8s^2\sqrt{2\frac{s-1}{s}}\log(2s-1)+o(1)\right)+\text{\rm O}(X\log^{-20}X),$$
where the implied constants do not depend on $D$.

We now compute a lower bound for the expression in (A.5) as in [Iw]. Since $D_1$ in our case is $1$ or $2$, the expression $\Omega_D$ in (4.8) of [Iw] becomes
$$
\Omega_D=c\cdot \sum_{\Sb{\delta}\\ {2|D\delta}\endSb }\sum_{\Sb{L\in\Cal L}\\{(\delta L + a,Q\delta)=1}\endSb}\frac{1}{\phi(Q\delta)}
\tag A.7
$$
where $c$ is a constant not depending on $D$ (coming from the products over $p|D_1$ in (4.8) of [Iw])
and $\delta= 1$ or $2$ as in Table 1.   Note that the innermost sum of the expression in
(4.8) is $\gg_{\epsilon} D^{-\epsilon}$ for $\epsilon>0$.  This follows from $|\Cal L|=\prod_{p|D_2}(p-1)/2\gg_{\epsilon'}D^{1-\epsilon'}$.  Define
$$
\tilde{\Omega}_D=\sum_{\Sb{\delta}\\ {|Q\delta|\leq
\log^{15}X}\\ {2|D\delta}\endSb}\sum_{\Sb{L\in\Cal L}\\{(\delta L + a,Q\delta)=1}\endSb}\frac{1}{\phi(Q\delta)}
$$
and note that, since $\delta\leq 2$ and $D\leq \log X$ in our case, 
$$
\aligned
|\Omega_a-\tilde{\Omega}_a|\leq\sum_{\Sb{\delta}\\ {Q\delta>\log^{15}X}\\ {p|\delta\Rightarrow p|D}\endSb}\frac{Q}{\phi(Q\delta)}
&<|8D|\cdot\sum_{Q\delta>\log^{15} X}\frac{1}{\sqrt{Q\delta}\sqrt{\phi(Q)}}\\
&<\frac{|8D|}{\log^{7.5} X}\\
&\leq \frac{1}{\log^{6} X}
\endaligned
$$
Combined with Theorem 1 of [Iw], this gives us the following bounds for $S_1(\phi, X, a)$ where $D\leq\log
X$ and $\delta=1$ or $2$:
$$
\aligned
S_1&\geq\theta\sqrt{\frac{2e^{\gamma}}{\pi}}C_0\cdot\tilde{\Omega}_a\frac{X}{(\log X)^{3/2}}(1+o(1))+\text{\rm O}(X\log^{-20}X)\\
&=\theta\Psi_D\Omega_D\frac{X}{(\log X)^{3/2}}(1+o(1))+\text{\rm O}(X\log^{-6}X) \\
\endaligned
$$
where the implied constants do not depend on $D$,
$$
\theta=\sup_{1<s<4/3}\left(\int_1^s\frac{dt}{\sqrt{t(t-1)}}-8s^2\sqrt{\frac{2(s-1)}{s}}\log(2s-1)\right),
$$
$$
C_0=\Psi_D=\sqrt{\frac{2}{\pi}}\prod_{p|2Da}\left(1-\frac{1}{p}\right)^{-1/2}\prod_{\Sb {p\not \; |2Da}\\{\left(\frac{k(-D)}{p}\right)=-1}\endSb}
\left(1-\frac{1}{(p-1)^2}\right)\prod_{p\not \; |2Da}\left(1-\frac{1}{p}\right)^{-\frac{1}{2}\left(\frac{k(-D)}{p}\right)}\gg_{\epsilon}D^{-\epsilon}$$
and $\Omega_D\gg_{\epsilon} D^{-\epsilon}$ as well for $\epsilon>0$.  This gives us the desired generalization of Iwaniec's theorem to Theorem~A.1.

\enddocument